\documentclass[11pt]{article}
\usepackage{amsthm}
\usepackage{latexsym}
\usepackage{amssymb}
\usepackage{amsmath}
\usepackage{mathrsfs}
\usepackage{array}
\usepackage{euscript}
\usepackage{graphicx}
\usepackage[T1]{fontenc}
\newcommand{\C}{\mathbb{C}}

\newcommand{\R}{\mathbb{R}}

\renewcommand{\S}{\mathbb{S}}

\newcommand{\Z}{\mathbb{Z}}

\renewcommand{\d}{\mathrm{d}}
\newcommand{\notD}{{{\cal D}\hskip-0.65em / \,}}
\newcommand{\D}{\mathrm{D}}

\newcommand{\dvol}{\mathrm{dVol}}

\newcommand{\thorn}{\mbox{\th}}

\newtheorem{definition}{Definition}[section]
\newtheorem{theorem}{Theorem}
\newtheorem{proposition}[definition]{Proposition}
\newtheorem{corollary}[definition]{Corollary}
\newtheorem{lemma}[definition]{Lemma}
\newtheorem{remark}[definition]{Remark}
\begin{document}\mbox{}

\vspace{0.25in}

\begin{center}
\huge{The characteristic Cauchy problem for Dirac fields on curved backgrounds}

\vspace{0.25in}

\large{Dietrich H\"AFNER\footnote{Université de Grenoble 1, Institut Fourier - UMR5582, \\ \noindent 100 rue des Maths, BP 74 38402 St Martin d'Heres, France. \\ \noindent
    Dietrich.Hafner@ujf-grenoble.fr} \&
    Jean-Philippe NICOLAS\footnote{Universit\'e de Brest, Laboratoire de Math\'ematiques de Brest, UMR 6205, \\
    6 avenue Victor Le Gorgeu, 29238 Brest cedex 3,
    France. \\ \noindent Jean-Philippe.Nicolas@univ-brest.fr}}
\end{center}

\begin{abstract}
On arbitrary spacetimes, we study the characteristic Cauchy problem for Dirac fields on a light-cone. We prove the existence and uniqueness of solutions in the future of the light-cone inside a geodesically convex neighbourhood of the vertex. This is done for data in $L^2$ and we give an explicit definition of the space of data on the light-cone producing a solution in $H^1$. The method is based on energy estimates following L. H\"ormander \cite{Ho}. The data for the characteristic Cauchy problem are only a half of the field, the other half is recovered from the characteristic data by integration of the constraints, consisting of the restriction of the Dirac equation to the cone. A precise analysis of the dynamics of light rays near the vertex of the cone is done in order to understand the integrability of the constraints~; for this, the Geroch-Held-Penrose formalism is used.
\end{abstract}

\tableofcontents

\section{Introduction}

The characteristic Cauchy problem, or Goursat problem, is a Cauchy problem for a hyperbolic equation, with data set on a characteristic hypersurface. The well-posedness depends on the geometry of the characteristic hypersurface. In the typical example of the scalar wave equation on $\R_t \times \R_x^3$, specifying data on the characteristic hyperplane $t=x_1$ leads to non-unique solutions, whereas for data on the forward light-cone of the origin $\{ t = \vert x \vert \}$, the solution exists and is unique in the future of the cone (but not in its past). In the best cases, the well-posedness will always be on one side of the hypersurface, its future or its past, unless we work on a spatially compact spacetime. A remarkable feature of the Goursat problem is that fewer data are necessary on a characteristic hypersurface than on a spacelike slice~; the remaining data can be recovered by integration of the restriction of the equation to the null hypersurface, referred to as {\it constraints}.

For the linear scalar wave equation on general globally hyperbolic curved spacetimes, the question of existence and uniqueness is well understood. A whole chapter of F.G. Friedlander's book \cite{Fri75} is devoted to an integral formulation of the solution for data on a light-cone using techniques due to Leray and Hadamard. Lars H\"ormander \cite{Ho} has proved global well-posedness for spatially compact spacetimes using a simple and natural method based on energy estimates (there is also a paper by J.-P. Nicolas \cite{Ni} extending H\"ormander's result to metrics of weak regularity). The much more delicate case of quasi-linear hyperbolic equations (including first order equations) was addressed by A. Rendall in \cite{Re} where he established local existence results with data on two intersecting null hypersurfaces and discussed applications to general relativity~; he also treated the case of data on a lightcone, for a special class of second order quasi-linear equations for which he could apply Friedlander's results to an associated linear equation. For spinorial zero rest-mass field equations, an integral Kirchhoff-d'Adh\'emard formula was obtained in the flat case by R. Penrose in 1963 \cite{Pe63} (see also \cite{PeRi} Vol. 1, Section 5.11) for data on a light-cone. In the curved case, Friedlander's approach had not been applied to spinorial equations until the very recent paper by J. Joudioux \cite{Jo} which completely generalizes Penrose's work to the general curved situation. H\"ormander's method was used in recent contributions by D. H\"afner \cite{Ha} to solve a Goursat problem for Dirac fields on the Kerr metric, with data on a characteristic surface generated by two congruences of outgoing and incoming null geodesics, and by L.J. Mason and J.-P. Nicolas \cite{MaNi} to construct scattering theories for Dirac, Maxwell and scalar fields via conformal methods.

The question of the regularity of the solutions and its control in terms of the regularity of the data is strikingly more difficult than for the ordinary Cauchy problem, particularly so when data are specified on a light-cone. This is due to the fact that null data, as was mentioned above, are only a part of the field, the remaining parts being obtained by integration of the equation restricted to the null hypersurface~; the control of the regularity of the solution in terms of the null data is consequently trickier, because already on the null hypersurface, the regularity of the full field depends in a complicated way on that of the null data. The additional difficulty on a light-cone is due to the singularity at the tip. Friedlander's book gives a condition ensuring smooth solutions for scalar waves and Rendall's work describes a way of controlling the Sobolev norm of the solution at any given order using the Whitney extension theorem (see for example H\"ormander \cite{HoBook}) but without trying to obtain the best space of null data for a given regularity of the solution. To our knowledge, a precise study of intermediate regularities in this sense is to this day missing. This work is a step in this direction.

We study the Goursat problem for the Dirac equation on a curved background, with data on a future light-cone. We work locally in a geodesically convex, globally hyperbolic neighbourhood of the vertex\footnote{A geodesically convex neighbourhood is also referred to as a normal convex neighbourhood.}, we therefore do not need to make any global hypothesis on our spacetime. We find the space of data on the cone for which the problem has a unique $H^1$ solution~; this is our main theorem (theorem \ref{T2}). Using density arguments, we infer a minimum regularity existence and uniqueness result given in theorem \ref{T1}. The general strategy of the proof is similar to that of H\"afner \cite{Ha} and uses the ideas developed by H\"ormander \cite{Ho} for the wave equation~: the Goursat problem is solved by constructing a trace operator on the null hypersurface, interpreting it by means of energy estimates both ways as a bounded linear operator that is one-to-one with closed range, then showing that it is onto. In particular, we prove a version for Dirac's equation of a general result in \cite{Ho}~: the well-posedness of the Cauchy problem on a rough (Lipschitz) hypersurface, on a spatially compact spacetime (this is given in theorem \ref{CauchyLipschitz}).

The main difficulty in the control of the $H^1$ regularity by an adequate space of characteristic data is the analysis of the operator which solves the constraint equations. This is due mostly to the singularity of the cone (this step is trivial in the case of two intersecting null hypersurfaces) and made formally more complicated by the spinorial nature of the equation. The constraints are typically considered as transport equations along the null generators of the cone with data at the vertex~; they are however singular at the tip of the cone. The proof of their integrability relies on a rescaling of a natural Newman-Penrose tetrad on the cone (which amounts in effect to defining a metric on the cone with its vertex blown up) and on a detailed study of the geometry of the null geodesic congruence on the cone. We show that the constraints have a unique solution that is bounded at the vertex and that it satisfies the direction-dependent matching conditions at the tip which characterize the continuity of the complete field. Since our equation is spinorial, we choose to perform our calculations using the Newman-Penrose formalism, or rather its compacted form~: the compacted spin-coefficient formalism (also referred to as the Geroch-Help-Penrose formalism, see \cite{PeRi} Vol.1 section 4.12). Our calculations using this formalism are similar to the detailed study of the geometry of null cones in Klainerman-Nicol\`o \cite{KlaNi}. Their approach is based on a double null foliation, which is similar in spirit to the Geroch-Held-Penrose formalism but without the spinorial aspect which is crucial to us.

The paper is organized as follows~:
\begin{itemize}
\item section \ref{GeBa} contains some useful geometrical background, a presentation of the Dirac equation and its conserved quantity as well as its description in the Newman-Penrose formalism~;
\item the main geometric objects and constructions of the paper are presented in section \ref{GeFra}, among which coordinate systems, Newman-Penrose tetrads, $3+1$ decomposition of the geometry, structure of the cone, blow-up of its tip and construction of a metric on the blown-up cone, function spaces on the cone~;
\item our three theorems are given in section \ref{Re} with first the $L^2$ result and an analogous result for an equation with source, the integrability of the constraints, the $H^1$ result and finally the Cauchy problem on a Lipschitz hypersurface~;
\item the proofs of the three main theorems are in section \ref{Pro}~;
\item the three appendices contain the proof of the integrability of the constraints and other fairly technical aspects of the paper~; appendix \ref{ConstraintEqSolve} presents the details of the proof of the integrability of the constraints with first a description in flat spacetime where the difficulties are already present and then the analysis of the geometry of null rays on the cone necessary for solving the constraints, appendix \ref{MoreSpinCoeffs} contains technical results in a similar spirit which are useful for showing that smooth functions are in the domain in $L^2$ of the operator solving the constraints, finally appendix \ref{CompactedNP} gives a short description of the compacted spin-coefficient formalism with some useful remarks on the Dirac equation.
\end{itemize}

{\bf Notations.} Many of our equations will be 
expressed using the two-component spinor notations and abstract index
formalism of R. Penrose and W. Rindler \cite{PeRi}. Abstract indices are
denoted by light face latin letters, capital for spinor indices and
lower case for tensor indices. Concrete indices defining components in reference to a basis are
represented by bold face latin letters. Concrete spinor indices,
denoted by bold face capital latin letters, take their values
in $\{ 0,1 \}$ while concrete tensor indices, denoted by bold face
lower case latin letters, take their values in $\{ 0,1,2,3
\}$. 

We will work with descriptions of the Dirac equation both in terms of Dirac spinors and Weyl (or half) spinors. For a complete account of the relations between Weyl and Dirac spinors, see \cite{PeRi} or \cite{Ni}.

\section{Geometrical background} \label{GeBa}

\subsection{Geodesic convexity, global hyperbolicity and spin structure} \label{GHSS}

In this work, we shall consider general Lorentzian manifolds and work locally in a geodesically convex and globally hyperbolic neighourhood of a point which always exists on a general smooth Lorentzian spacetime (see Friedlander \cite{Fri}, theorem 4.4.1 p. 147).
\begin{definition} \label{GeodConvex}
A domain is said to be geodesically convex if between any two points of the domain there exists a unique geodesic that is entirely contained in the domain.
\end{definition}
\begin{remark}
Such domains are called convex normal in the classic textbook by Wald \cite{Wa}, whereas Friedlander in \cite{Fri75} calls them geodesically convex.
\end{remark}
We recall here the definition of global hyperbolicity and some of its important consequences in dimension $4$, particularly regarding spinors.

A globally hyperbolic spacetime is a pair $({\cal M},g)$ where (see Geroch \cite{Ge3} for more details)~:
\begin{itemize}
\item $\cal M$ is a real $4$-dimensional, smooth, oriented, time-oriented manifold~;
\item $g$ is a smooth metric on $\cal M$ of Lorentzian signature $+\, - \, - \, -$~;
\item there exists a global time function $t$ on $\cal M$ such that the level hypersurfaces ${\Sigma}_t$ of $t$ are Cauchy hypersurfaces.
\end{itemize}
The time function $t$ may be in addition assumed smooth (see Bernal-Sanchez \cite{BeSa}). Recall that a smooth time function is a smooth scalar function $t$ on $\cal M$ such that $\nabla^a t$ is a future-oriented timelike vector field over $\cal M$~; here $\nabla$ denotes the Levi-Civita connection on $({\cal M},g)$.

Global hyperbolicity has at least two important consequences in $4$ dimensions (see \cite{BeSa, Ge1, Ge2, Stie}). First, the level hypersurfaces $\Sigma_t$ of the time fonction $t$ are all diffeomorphic to a given smooth $3$-surface $\Sigma ={\Sigma}_0$ via the flow of the vector field $\nabla^a t$. Second, $\cal M$ admits a spin-structure. We denote by $\mathbb{S}_A$ and $\bar{\mathbb{S}}_{A'}$ the bundles of left and right spinors on $\cal M$. These are $2$-component spinors, or Weyl spinors. The Weyl spinor bundles are endowed with symplectic forms $\varepsilon_{AB}$ and $\varepsilon_{A'B'}$ which are conjugates of one another. They are used to raise and lower spinor indices (meaning that they provide isomorphisms between the spin-bundles $\mathbb{S}^A$ and $\bar{\mathbb{S}}^{A'}$ and their duals $\mathbb{S}_A$ and $\bar{\mathbb{S}}_{A'}$). The bundle of Dirac spinors is defined as
\[ \mathbb{S}_\mathit{Dirac} := \mathbb{S}_A \oplus \bar{\mathbb{S}}^{A'} \, .\]
It is equipped with an $SL(2,\C )$
invariant inner product expressed as
\begin{equation} \label{SympFormSDirac}
\left( \Psi \, ,~ \Xi \right) := i \bar\rho_{A'} \chi^{A'} -i \phi_{A} \bar\eta^{A} \, ,~\mbox{ where } \Psi = \phi_A \oplus \chi^{A'} \mbox{ and } \Xi = \rho_A \oplus \eta^{A'} \, .
\end{equation}
\begin{remark}
The Clifford product by a vector is skew for the symplectic product on Dirac spinors, but for the sesqui-linear $2$-form \eqref{SympFormSDirac}, which is the symplectic form applied to $\Psi$ and the dual of the conjugate of $\Xi$, the Clifford product by any real vector is symmetric.
\end{remark}

The tangent bundle to $\cal M$ and the metric can be recovered from the Weyl-spinor bundles and the $\varepsilon$ symplectic forms~:
\[ T {\cal M} \otimes \C = \mathbb{S}^A \otimes \bar{\mathbb{S}}^{A'} \, ,\]
(a rigorous abstract index notation should in fact be $T^a {\cal M} \otimes \C = \mathbb{S}^A \otimes \bar{\mathbb{S}}^{A'}$, the vector index $a$ corresponding to the two spinor indices $A$ and $A'$ clumped together), the real tangent bundle consists of the hermitian part of $\mathbb{S}^A \otimes \bar{\mathbb{S}}^{A'}$ and
\[ g_{ab} = \varepsilon_{AB} \varepsilon_{A'B'} \, .\]
We can perform a $3+1$ decomposition of the geometry based on the time function $t$. We normalize the gradient of $t$ so that its square norm equals $2$ (instead of a more usual $1$, this is for later convenience in the expression of the hermitian norm of spinors defined using this vector field)
\[ {\cal T}^a := \sqrt{\frac{2}{g(\nabla t , \nabla t )}} \nabla^a t \, . \]
The metric $g$ can then be decomposed as follows~:
\[ g = \frac{N^2}2 \d t^2 - h(t) \]
where $-h(t)$ is the metric induced by $g$ on $\Sigma_t$ and the lapse function $N$ is defined by
\[ {\cal T}_a \d x^a = N \d t \, , ~\mbox{or equivalently}~g(\nabla t , \nabla t ) = \frac{2}{N^2} \, .\]
Note that this decomposition, and more particularly the choice of product structure ${\cal M} = \R_t \times \Sigma$ associated with the integral curves of ${\cal T}$, fixes the meaning of the vector $\partial / \partial t$ as
\[ \frac{\partial}{\partial t} = \frac{N}{2} {\cal T} = \frac{\nabla t}{g( \nabla t \, , \nabla t )} \, .\]
\begin{definition}
The timelike vector ${\cal T}$ endows the bundle of Dirac spinors with a positive definite hermitian product~:
\begin{equation} \label{DiracSpinorInnerProduct}
\langle \Psi \, ,~ \Xi \rangle := \frac{1}{\sqrt{2}} \left( {\cal T} . \Psi \, ,~ \Xi \right) = {\cal T}^{AA'} \phi_A \bar{\rho}_{A'} + {\cal T}_{AA'} \chi^{A'} \bar{\eta}^A \, ,~ \Psi = \phi_A \oplus \chi^{A'} \, ,~\Xi = \rho_A \oplus \eta^{A'} \, .
\end{equation}
In particular we denote
\begin{equation} \label{DiracSpinorNorm}
\left| \Psi \right| = \langle \Psi , \Psi \rangle^{1/2} \, .
\end{equation}
\end{definition}
\begin{remark}
The Clifford product by a real vector co-linear to $\cal T$ (resp. orthogonal to $\cal T$) is self-adjoint (resp. skew-adjoint) for the inner product \eqref{DiracSpinorInnerProduct}.
\end{remark}
\begin{remark}
Note that this definition can be naturally restricted to each of the bundles $\mathbb{S}_A$ and $\bar{\mathbb{S}}^{A'}$ and extended to $\mathbb{S}^A$ and $\bar{\mathbb{S}}_{A'}$ .
\end{remark}

\subsection{Dirac's equation, conserved quantity}

We consider the charged Dirac equation associated with an electromagnetic vector-potential $\Phi^a$, for a particle of mass $m$ and charge $q$
\begin{equation} \label{DirEq}
\left. \begin{array}{ccl}{\big( \nabla^{AA'} -iq \Phi^{AA'} \big) \phi_A} &  = & {\frac{m}{\sqrt{2}} \chi^{A'} ,} \\
{\big( \nabla_{AA'} -iq \Phi_{AA'} \big) \chi^{A'}} &  = & -{\frac{m}{\sqrt{2}} \phi_{A} \, .} \end{array} \right\}
\end{equation}
This is an expression of Dirac's equation in terms of Weyl spinors, it has the form of two charged Weyl equations, one for helicity $1/2$ and the other for helicity $-1/2$, coupled by the mass. It is usual to understand Dirac's equation as an equation on the Dirac spinor $\Psi = \phi_A \oplus \chi^{A'}$. Contrary to the more elementary spin-bundles $\mathbb{S}_A$ and $\bar{\mathbb{S}}^{A'}$, the Dirac spinor bundle $\mathbb{S}_\mathit{Dirac}$ is stable under the action of the Clifford product by vectors and of the Dirac operator (see for example \cite{Ni} or \cite{PeRi} for a description of the relations between the Dirac operator acting on Weyl spinors and on Dirac spinors)
\[ \notD \, :~ \mathbb{S}_A \oplus \bar{\mathbb{S}}^{A'} \longrightarrow \mathbb{S}_A \oplus \bar{\mathbb{S}}^{A'} \, ,~ \notD = \left( \begin{array}{cc} 0 & {i \sqrt{2} \nabla_{AA'} } \\ {-i\sqrt{2} \nabla^{AA'}} & 0 \end{array} \right) \, . \]

Equation \eqref{DirEq} has a conserved current given by the future-oriented causal vector (see \cite{PeRi})
\[ J^a = \phi^A \bar{\phi}^{A'} + \bar{\chi}^A \chi^{A'} \, .\]
For a given spacelike or characteristic hypersurface $\cal S$, the flux of $J^a$ through $\cal S$ is the integral over $\cal S$ of the Hodge dual $J_a \d^3 x^a$ of the $1$-form $J_a \d x^a$~:
\begin{gather}
{\cal E}_{\cal S} := \int_{\cal S} J_a \d^3 x^a , \,J_a \d^3 x^a := * J_a \d x^a . \label{Flux}
\end{gather}
The flux \eqref{Flux} can be understood as
\begin{equation} \label{FluxAlternative}
{\cal E}_{\cal S} = \int_{\cal S} J_a V^a \d \sigma_{\cal S}
\end{equation}
where $V^a$ is orthogonal to ${\cal S}$, $d\sigma_{\cal S}=W\lrcorner \dvol^4$, $\dvol^4$ being the $4$-volume measure associated with the metric $g$ and $W^a$ a transverse vector field to ${\cal S}$ such that $V_a W^a=1$. If ${\cal S}$ is spacelike we can take $V^a=W^a=\nu^a$, the unit future oriented normal vector field to ${\cal S}$. If ${\cal S}$ is characteristic, the vector $V^a$ is a null vector field normal and tangent to $\cal S$ and we can chose $W^a$ to be a null vector field transverse to $\cal S$, this provides the beginning of the construction of a Newman-Penrose tetrad (see Section \ref{NewPenform}).

Provided we have a local orthonormal frame $\{ e_0 \, ,~ e_1 \, ,~ e_2 \, , ~ e_3 \}$, the Dirac operator $\notD$ can be decomposed as follows
\begin{equation} \label{3+1Dirac1}
\notD = \sum_{\alpha =0}^3 e_\alpha . \nabla_{e_\alpha}
\end{equation}
and equation \eqref{DirEq} reads, 
\begin{equation} \label{3+1Dirac2}
( \notD + {\cal P} ) \Psi =0 \, ,
\end{equation}
where $\cal P$ is a smooth potential involving the mass and charge terms.

\subsection{Newman-Penrose formalism}
\label{NewPenform}

We will make an essential use in this work of the expression of equation \eqref{DirEq} in the Newman-Penrose formalism. This formalism is based on the choice of a null tetrad, i.e. a set of four vector fields $l^a$, $n^a$, $m^a$ and $\bar{m}^a$, the first two being real and future oriented, $\bar{m}^a$ being the complex conjugate of $m^a$, such that all four vector fields are null and $m^a$ is orthogonal to $l^a$ and $n^a$, that is to say
\[ l_a l^a = n_a n^a = m_a m^a = l_a m^a = n_a m^a = 0 \, .\]
The tetrad is said to be normalized if in addition
\[ l_a n^a = 1 \, ,~ m_a \bar{m}^a =-1 \, .\]
Such a tetrad induces a decomposition of the metric as follows (see \cite{PeRi} Vol.1 p.119)~:
\[ g_{ab} = l_a n_b + n_a l_b - m_a \bar{m}_b - \bar{m}_a m_b \, .\]
To  a given Newman-Penrose tetrad we can associate a spin-frame $\{o^A ,\iota^A \}$, i.e. a local basis of the spin-bundle $\mathbb{S}^A$, defined uniquely up to an overall sign factor by
\[ o^A \bar{o}^{A'} = l^a \, ,~\iota^A \bar{\iota}^{A'} = n^a \, ,~o^A
\bar{\iota}^{A'} = m^a \, ,~\iota^A \bar{o}^{A'} = \bar{m}^a \, .\]
If the Newman-Penrose tetrad is normalized, then the corresponding spin-frame satisfies $o_A \iota^A =1$ and is said to be unitary.

Let $\phi_0$ and $\phi_1$ be the components of $\phi_A$ in
$\{o^A , \iota^A \}$, and $\chi_{0'}$ and $\chi_{1'}$ the components of
$\chi_{A'}$ in $(\bar{o}^{A'} , \bar{\iota}^{A'} )$~:
\[ \phi_0 = \phi_A o^A~,~~\phi_1 = \phi_A \iota^A ~,~~ \chi_{0'} =
\chi_{A'} \bar{o}^{A'} ~,~~ \chi_{1'} = \chi_{A'} \bar{\iota}^{A'} \,
.\]
The Dirac equation takes the form (see for example \cite{Ch})
\begin{eqnarray}
\label{NewmanDirac}
 \left. \begin{array}{l}
{ n^\mathbf{a} (\partial_\mathbf{a}-iq\Phi_{\bf{a}}) \, \phi_0 - m^\mathbf{a}
(\partial_\mathbf{a}-iq\Phi_{\bf{a}})  \, \phi_1 + (\mu - \gamma )\phi_0 + (\tau - \beta )
\phi_1 = \frac{m}{\sqrt{2}} \chi_{1'} \, , } \\ \\
{ l^\mathbf{a} (\partial_\mathbf{a}-iq\Phi_{\bf{a}})  \, \phi_1 - \bar{m}^\mathbf{a}
(\partial_\mathbf{a}-iq\Phi_{\bf{a}})  \, \phi_0 + (\alpha - \pi )\phi_0 + (\varepsilon -
\rho ) \phi_1 = - \frac{m}{\sqrt{2}} \chi_{0'} \, , } \\ \\
{ n^\mathbf{a} (\partial_\mathbf{a}-iq\Phi_{\bf{a}})  \, \chi_{0'} - \bar{m}^\mathbf{a}
(\partial_\mathbf{a}-iq\Phi_{\bf{a}})  \, \chi_{1'} + (\bar{\mu} - \bar{\gamma} )\chi_{0'}
+ (\bar{\tau} - \bar{\beta} ) \chi_{1'} = \frac{m}{\sqrt{2}} \phi_{1} \,
, } \\ \\
{ l^\mathbf{a} (\partial_\mathbf{a}-iq\Phi_{\bf{a}})  \, \chi_{1'} - m^\mathbf{a}
(\partial_\mathbf{a}-iq\Phi_{\bf{a}})  \, \chi_{0'} + (\bar{\alpha} - \bar{\pi} )\chi_{0'}
+ (\bar{\varepsilon} - \bar{\rho} ) \chi_{1'} = - \frac{m}{\sqrt{2}}
\phi_{0} \, .} \end{array} \right\}
\end{eqnarray}
The $\mu,\gamma$ etc. are the spin coefficients which are decompositions of the connection based on the vectors of the null tetrad, for instance, $\mu=-\bar{m}^a\delta n_a$, where $\delta=m^a\nabla_a$. For the formulae of the spin coefficients and details about the Newman-Penrose formalism see \cite{PeRi}.

A Newman-Penrose tetrad $(l,n,m,\bar{m})$ is said to be adapted to the foliation $\{ \Sigma_t \}_t$ if it satisfies $l^a + n^a = {\cal T}^a$. The advantage of a tetrad adapted to the foliation is that the expression of the hermitian product on Dirac spinors becomes extremely simple. Let $\Psi = \phi_A \oplus \chi^{A'}$ and $\Xi = \rho_A \oplus \eta^{A'}$, denote the four components of $\Psi$ in the spin-frame $\{ o^A , \iota^A \}$ by
\begin{equation} \label{DiracSpinorComponents}
(\Psi_1, \Psi_2 , \Psi_3 , \Psi_4)=(\phi_0,\phi_1,\chi^{0'},\chi^{1'})=(\phi_0,\phi_1,\chi_{1'},-\chi_{0'}) \, ,
\end{equation}
then, since we have ${\cal T}^a = l^a + n^a = o^A \bar{o}^{A'} + \iota^A \bar{\iota}^{A'}$,
\begin{equation} \label{AdaptedInnerProduct}
\left| \Psi \right|^2 = \left| \Psi_1 \right|^2 + \left| \Psi_2 \right|^2 + \left| \Psi_3 \right|^2 + \left| \Psi_4 \right|^2
\end{equation}
and with analogous notations for $\Xi$~:
\[ \langle \Psi , \Xi \rangle = \Psi_1 \bar{\Xi}_1 + \Psi_2 \bar{\Xi}_2 + \Psi_3 \bar{\Xi}_3 + \Psi_4 \bar{\Xi}_4 \, .\]

\section{Geometrical framework} \label{GeFra}

On a smooth Lorentzian spacetime (${\cal M} \, ,~ g$), we consider a point $p_0$ and we work in a geodesically convex (see definition \ref{GeodConvex}) and globally hyperbolic neighbourhood $\Omega$ of $p_0$. Every point in a smooth spacetime admits such a neighbourhood (see Friedlander \cite{Fri}, theorem 4.4.1 p. 147).

\subsection{Important sets, optical functions and time function}

For a point $p \in \Omega$ the future light-cone of $p$ is defined as the set of points which are null separated from $p$ and in the future of $p$
\[ {\cal C}^+ (p) := \left\{ q \in \Omega \, ;~ \mbox{there exists a future-oriented null geodesic from } p \mbox{ to } q \right\} \cup \{ p \} \, , \]
we denote by ${\cal I}^+(p)$ the future chronological set of $p$ in $\Omega$
\[ {\cal I}^+ (p):= \left\{ q \in \Omega \, , q \neq p \, ;~ \mbox{there exists a future-oriented timelike geodesic from } p \mbox{ to } q \right\}\]
and by ${\cal J}^+ (p)$ the future causal set of $p$ in $\Omega$
\[ {\cal J}^+ (p) = {\cal I}^+ (p) \cup {\cal C}^+ (p) \, .\]
We define the analogous sets in the past in the natural way~: ${\cal C}^- (p)$, ${\cal I}^- (p)$ and ${\cal J}^- (p)$\footnote{We could also have defined these sets using curves which are not geodesics. In our case, both definitions coincide, see Wald \cite{Wa}, theorem 8.1.2, p. 191.}.

For the resolution of the Goursat problem, we work in a neighbourhood of the vertex of the cone ${\cal C}^+ (p_0)$ within ${\cal J}^+ (p_0)$, which does not need to be small, it can be large if $\Omega$ is itself large. We define this neighbourhood in two steps~: first we construct two optical functions and a time function on $\Omega$, then we use this time function to define the domain in which we shall solve the Goursat problem.

We construct two foliations of $\Omega$ by light-cones which we use to define optical functions.
\begin{definition}
Let $\zeta$ be a timelike geodesic passing through $p_0$, inextendible in $\Omega$. We define a time function $t$ on $\zeta$ as the metric length of $\zeta$ between $p_0$ and the points on $\zeta$, with a positive (resp. negative) sign if we are in the future (resp. past) of $p_0$. We put $T_1 = \mathrm{inf} \{ t \}$ and $T_2= \mathrm{sup} \{ t \}$ on $\zeta$ in $\Omega$. We denote for $t\in ]T_1,T_2[$,
\[ {\cal C}^\pm_t := {\cal C}^\pm (\zeta (t))  \mbox{ in } \Omega \, .\]
We consider the open, geodesically convex, globally hyperbolic subset of $\Omega$~:
\[ D := \bigcup_{T_1 <t_1 <t_2< T_2 } \left( {\cal J}^+(\zeta (t_1)) \cap {\cal J}^- (\zeta (t_2)) \right) \, .\]
We define two optical functions (i.e. solutions of the eikonal equation $g(\nabla \xi , \nabla \xi ) =0$) $u$ and $v$ on $D$ as follows~:
\begin{itemize}
\item $u:=t$ on ${\cal C}^+_t$~;
\item $v:=t$ on ${\cal C}^-_t$.
\end{itemize}
Then we extend the definition of the time function $t$ to the whole of $D$ by putting~:
\[ t := (u+v)/2 \, .\]
\end{definition}
\begin{remark}
Note that the coordinate $t$ as we have defined it on $D$ agrees with its initial definition on $\zeta$.
\end{remark}
\begin{remark}
The functions $u$ and $v$ are indeed solutions of the eikonal equation since they are constant on the future and past null cones respectively, their gradients are therefore othogonal to these null hypersurfaces whose orthogonal subspace is the span of a null vector.
\end{remark}
The function $t$ has the following property~:
\begin{proposition}
The function $t $ is a smooth time function on $D$.
\end{proposition}
{\bf Proof.} To prove the smoothness of $t$ on $D$, it is sufficient by smooth angular dependence to prove it on a well chosen $2$-surface $\cal S$ inside $D$. We define the surface $\cal S$ as follows. Consider a future oriented null geodesic $\gamma_0$ passing through $p_0$ with some affine parameter $s$. We parallel-transport the tangent vector $\gamma_0'(0)$ to $\gamma_0$ at $p_0$ along the geodesic $\zeta$. We obtain at each point $\zeta (t)$ a future oriented null vector and we consider the associated null geodesic $\gamma_{t}$ passing through $\zeta (t)$. The family of geodesics $\gamma_t$ spans within $D$ a smooth $2$ surface $\cal S$. When we rotate the direction of the initial null geodesic $\gamma_0$ through all the possible directions $\omega$ (quotiented by the antipodal relation~: $\omega \mapsto -\omega$), the resulting $2$-surface will change smoothly and define a smooth foliation of $D \setminus \zeta$ (thanks to the geodesic convexity of $\Omega$).

Now for a given geodesic $\gamma_0$, we work on the corresponding $2$-surface $\cal S$. First, note that for any $t$, the intersection $\left( {\cal C}^+_t \cup {\cal C}^-_t \right) \cap {\cal S}$ is the union of two smooth curves, one being $\gamma_t$. We denote the other by $\beta_t$. The geodesic $\zeta$ splits $\cal S$ into two halves~: a left and a right part. We change the definitions of $u$ and $v$ as follows. On the right part, we keep $u$ and $v$ as they are and on the left part, we exchange the roles of $u$ and $v$. Then the $u=t$ curves on $\cal S$ are the $\gamma_t$'s and the $v=t$ curves are the $\beta_t$'s. The essential remark is that if now we put $t := (u+v)/2$, this does not change the definition of $t$, which shows immediately that $t$ is smooth on $\cal S$ and thus on $D$.

Moreover the gradient of $t$ is $\nabla t = \frac12 (\nabla u + \nabla v)$ which is the sum of two future oriented null vectors on $D \setminus \zeta$ that are nowhere colinear and do not become so as one approaches $\zeta$~; it is therefore future oriented and timelike on $D \setminus \zeta$ and thus on $D$ by continuity. \qed
\begin{remark}
Note that in general, the level hypersurfaces of $t$ are not orthogonal to $\zeta$, i.e. $\zeta$ is not an integral curve of $\nabla t$.
\end{remark}

We perform the $3+1$ decomposition of the metric $g$ in $D$ based on the time function $t$ and we adopt the notations of section \ref{GHSS} for the elements of this decomposition. We denote by $\Sigma_t$ the level hypersurfaces of the function $t$ in $D$.

{\bf A technical choice.} Rescaling the variables $u$ and $v$ globally (and therefore $t$ as well) by a constant factor $k >0$ if necessary~:
\[ u \mapsto k u \, ,~ v \mapsto k v \, ,~ t \mapsto k t \, ,\]
we assume from now on that
\begin{equation} \label{LapseP0}
N (p_0 ) = \sqrt2 \, .
\end{equation}
This will provide important simplifications in the proof of the integrability of the constraint equations along the cone (see appendix \ref{ConstraintEqSolve}).
\begin{definition}
Let $0<T<T_2/2$, we shall henceforth work in the domain
\[ D_T := {\cal J}^+ (p_0) \cap \{ 0 \leq t \leq T \} \, .\]
\end{definition}
\begin{remark}
{\bf From now on, all the sets we shall consider will be restricted to $D_T$.} In particular, we keep the same notations for the sets defined above, such as ${\cal C}_t^+$ for instance, but we now consider only the restriction of these sets to $D_T$.
\end{remark}
In the next subsection, we define a local frame using the construction of our optical functions in $D$. This local frame will be singular at the curve $\zeta$~: more precisely, the frame vectors will be smooth on $\D_T \setminus \zeta$ and have direction dependent limits on $\zeta$. These limits and their relations to one another will be fundamental for our constructions. A natural way to deal with them is to blow up the curve $\zeta$ as the cylinder $[0,T]_t \times S^2_\omega$. This is done in the usual way~:
\begin{definition}
The domain $D_T$ with the curve $\zeta$ blown-up is defined as the set of pairs $(p,\gamma )$ where $p \in D_T$ and $\gamma$ is a future null geodesic passing through $p$ and through $\zeta$. When $p \in D_T \setminus \zeta$, there is only one such geodesic and the pair $(p , \gamma )$ can be unambiguously identified with $p$. When $p \in \zeta$ however, we have a whole $2$-sphere of possible geodesics $\gamma$. So this blow-up amounts to replacing each point $p$ on $\zeta$ by the sphere of future null geodesics from $p$. We denote by $\mathfrak{D}_T$ the domain $D_T$ with $\zeta$ blown up.

A similar construction can be performed on ${\cal C}^+_0$. The blow-up only affects the point $p_0$ which is replaced by the $2$-sphere of future null geodesics from $p_0$. We denote by $\mathfrak{C}$ the cone ${\cal C}^+_0$ with its vertex blown-up. It has the topology of a cylinder.
\end{definition}
We define a coordinate system in $\D_T$ based on the future lightcones.
\begin{definition}[Coordinate systems on $D_T$ and ${\cal C}_0^+$]
At $p_0$ we choose a smooth parametrization of future null geodesics by $\omega \in S^2$ and we extend it to all future null directions at points of $\zeta$ by imposing that $\omega$ remains constant when we transport a future null vector parallel along $\zeta$. Then we also impose that $\omega$ be constant along future null geodesics. Thus, future null geodesics emerging from a point $p \in \zeta$ are characterized by $t(p) \in [0,T]$ and $\omega \in S^2$ and a point $q$ on such a geodesic is characterised by $t(p)$, $\omega$ and $t(q)$. We choose to describe such a point $q$ by the coordinates $t:=t(q)$, $r:=t(q)-t(p)$ and $\omega$. The relation between $r$ and the coordinates $u$ and $v$ is simply
\[ r= t-u= (v - u)/2 \]
and thus
\[ u=t-r \, ,~ v=t+r \, .\]
On $\D_T$, $t$ varies in $[0,T]$, $r$ takes its values in $[0,t]$ for each $t$ and $\omega$ varies on $S^2$ and can be understood as parametrizing the $2$-surfaces orthogonal to $\nabla u$ and $\nabla v$ (the intersections of the future and past null cones emanating from $\zeta$).

When working on ${\cal C}_0^+$, we shall use the coordinates $v$ and $\omega$ as well as $t=r$ and $\omega$. Of course on ${\cal C}_0^+$, we have $v = t+r = 2t = 2r$.
\end{definition}
\begin{remark} \label{OmegaConstantGradu}
By lemma \ref{NablauGeod} and the fact that the orthogonal subspaces to the outgoing cones contain only one null direction, the future null geodesics on the cones ${\cal C}^+_t$ are the integral lines of $\nabla u$. In particular, $\omega$ is constant along the integral lines of $\nabla u$.
\end{remark}
We define a notion of smoothness on $\mathfrak{C}$ and $\mathfrak{D}_T$.
\begin{definition}[special regularity] \label{SmoothOnBlownUpCone}
A scalar function on $\mathfrak{D}_T$ is said to be s-smooth or s-regular (for specially smooth or regular) if it is ${\cal C}^\infty$ for the differential structure on
\[ \mathfrak{D}_T = \{ (u,v) \, ;~ 0 \leq u+v \leq 2T \, ,~ 0\leq u \leq v \} \times S^2 \, .\]
A vector field on $\mathfrak{D}_T$ is said to be s-smooth if, when contracted with any \underline{smooth $1$-form on $D_T$}, it gives an s-smooth function on $\mathfrak{D}_T$. The s-smooth $1$-forms on $\mathfrak{D}_T$ are those which, when contracted with s-smooth vector fields, give s-smooth functions. The definition of s-smooth tensor fields follows naturally. We define analogously s-smooth spinor fields~: a spinor field is s-smooth if its contraction with a smooth dual spinor field on $D_T$ gives an s-smooth function.

The s-smooth objects on $\mathfrak{C}$ are the restrictions to $\mathfrak{C}$ of s-smooth objects on $\mathfrak{D}_T$.
\end{definition}
\begin{remark} \label{RegRem}

\noindent
\begin{enumerate}
\item  The definition of s-regular vector fields is equivalent to the following property~: when applied to \underline{smooth functions on $D_T$} (and not to s-smooth functions on $D_T$!) the vector field produces an s-smooth function on $\mathfrak{D}_T$ (this amounts to contracting it with exact smooth $1$-forms~; we can always choose a basis of such $1$-forms). 
\item The notion of s-regularity on $\mathfrak{C}$ is not restricted to tangent objects. A vector field transverse to $\mathfrak{C}$ can be s-smooth on $\mathfrak{C}$.
\item Any object that is smooth on $D_T$ is s-smooth on $\mathfrak{D}_T$ and on $\mathfrak{C}$.
\item In general if we work with a Newman-Penrose tetrad which is adapted to the cone, not all frame vectors will be regular on the cone or the blown up cone in the usual sense, but they will be s-regular (see Section \ref{NPT})~; a typical example is that of the vectors $m$ and $\bar{m}$ in the Newman-Penrose tetrad we shall define in the next section.
\item Smoothness on $\mathfrak{D}_T$ (resp. $\mathfrak{C}$) and s-smoothness are distinct notions. The example above shows that s-smoothness does not imply smoothness on the blown-up domains. Conversely, if a scalar function $f$ is smooth on $\mathfrak{D}_T$ (resp. $\mathfrak{C}$) then so is its differential, but $\d f$ is usually not s-smooth.
\end{enumerate}
\end{remark}

\subsection{Newman-Penrose tetrad}
\label{NPT}

We now construct on $D_T$ a Newman-Penrose tetrad which is adapted to the foliation. The vector ${\cal T}^a$, future oriented, normal to the foliation $\{ \Sigma_t \}_{0\leq t \leq T}$ of $D_T$ is given by
\[ {\cal T}^a \partial_a = \sqrt{\frac{2}{g(\nabla t , \nabla t)}} \nabla t = N \nabla t \, .\]
Note that
\[ g(\nabla t , \nabla t) = \frac12 g(\nabla u , \nabla v) \, .\]
\begin{definition}
We define at each point of $D_T \setminus \zeta$ the vectors $l$ and $n$ of our Newman-Penrose tetrad by
\[ l =  \frac{\nabla u}{\sqrt{2 \, g(\nabla t , \nabla t)}} = \frac{N}{2} \nabla u \, ,~n = \frac{\nabla v}{\sqrt{2 \, g(\nabla t , \nabla t)}} = \frac{N}{2} \nabla v \, .\]
This defines $l$ and $n$ as s-smooth future null vector fields on $\mathfrak{D}_T$, satisfying $g(l,n)=1$, $l$ being tangent (and orthogonal) to the cones ${\cal C}^+_t$ and $n$ to the cones ${\cal C}^-_t$. A choice of $l$ and $n$ fixes the vector $m$ uniquely up to a complex factor of modulus $1$. We make a choice of $m$ so that it extends as an s-smooth vector field on $\mathfrak{D}_T$ (see remark \ref{TwoCharts}). The vectors $m$ and $\bar{m}$ are by construction tangent to the $2$-surfaces where both $u$ and $v$ (equivalently both $t$ and $r$) are constant, i.e. the intersections ${\cal C}^+_u \cap {\cal C}^-_v$.

We then extend this Newman-Penrose tetrad as an s-smooth Newman-Penrose tetrad on $\mathfrak{D}_T$~: we define a Newman-Penrose tetrad at each point $\zeta (t)$ for each direction $\omega$ by taking the limit of $(l,n,m,\bar{m} )$ along the integral curve of $\nabla u$ (null geodesic) on ${\cal C}_t^+$ corresponding to the direction $\omega$.

We consider at each point of $\mathfrak{D}_T$ the unique (modulo overall sign) spin-frame $\{ o^A \, ,~\iota^A \}$ associated to our null tetrad at this point~; we choose the sign consistently so that $o^A$ and $\iota^A$ are smooth over $\mathfrak{D}_T$.
\end{definition}
\begin{remark} \label{TwoCharts}
Note that the vectors $m$ and $\bar{m}$ cannot in fact be s-smooth on the whole $\mathfrak{D}_T$ since they are tangent to the surfaces of constant $u$ and $v$ and are normalized by $m_a \bar{m}^a =-1$. What must be done is to work with two charts on $\mathfrak{D}_T$ and in each chart a global choice of $m$ and $\bar{m}$ can be done. This is the usual problem linked with spherical coordinates (i.e. with zeros of vector fields tangent to a $2$-sphere) and we shall systematically ignore it.
\end{remark}
\begin{remark}
It is important to note that by continuity, the tetrads defined on $\zeta$ are still adapted to the foliation, i.e. they satisfy $l^a + n^a = {\cal T}^a$.
\end{remark}
The notion of s-regularity entails a useful property involving the vector $l$~:
\begin{lemma} \label{Sregl}
\begin{enumerate}
\item The vector field $l$ is smooth in the usual sense on $\mathfrak{D}_T$ and on $\mathfrak{C}$ (similarly $n$ is smooth on $\mathfrak{D}_T$ but not on $\mathfrak{C}$ since it is not tangent to the cone).
\item If a tensor or spinor field $T$ is s-smooth on $\mathfrak{D}_T$ (resp. $\mathfrak{C}$), then so is $\nabla_l T$.
\end{enumerate}
\end{lemma}
{\bf Proof.} \begin{enumerate}
\item The functions $u$ and $v$ are smooth on $\mathfrak{D}_T$ and therefore so are $\d u$ and $\d v$. This and the fact that $l$ is tangent to $\mathfrak{C}$ give the first property.
\item We first prove the result for a scalar field. Let $f$ an s-smooth function on $\mathfrak{D}_T$ (resp. $\mathfrak{C}$), this means that $f$ is smooth in the usual sense on $\mathfrak{D}_T$ (resp. $\mathfrak{C}$) and so is $\d f$. The first part of the lemma gives the result. Now consider $V$ an s-smooth vector field on $\mathfrak{D}_T$ (resp. $\mathfrak{C}$)~; all we need to show is that for any smooth $1$-form $\omega$ on $D_T$, $\omega_a \nabla_l V^a$ is s-smooth on $\mathfrak{D}_T$ (resp. $\mathfrak{C}$). Let $\omega$ a smooth $1$-form on $D_T$, we have
\begin{eqnarray*}
\omega_a \nabla_l V^a &=& \nabla_l (\omega_a V^a ) - V^a \nabla_l \omega_a \\
&=& \nabla_l (\omega_a V^a ) - V^a l^b \nabla_b \omega_a \, .
\end{eqnarray*}
The function $\omega_a V^a$ is s-smooth on $\mathfrak{D}_T$ (resp. $\mathfrak{C}$) and therefore, by what we just proved, so is $\nabla_l (\omega_a V^a )$. Now the $2$-form $\nabla_b \omega_a$ is smooth on $D_T$ which entails that $V^a l^b \nabla_b \omega_a$ is an s-smooth function on $\mathfrak{D}_T$ (resp. $\mathfrak{C}$). This proves the s-regularity of $\omega_a \nabla_l V^a$. The result extends to tensor fields by the Leibnitz rule and to spinor fields by an analogous argument. \qed
\end{enumerate}
\begin{remark}
This is in essence saying that if an object is s-smooth on $\mathfrak{C}$ (resp. $\mathfrak{D}_T$) then it is smooth on every future null geodesic $\gamma$ emanating from $p_0$ (resp. $\zeta$) right down to $p_0$ (resp. to the origin of $\gamma$ on $\zeta$) and its restriction to $\gamma$ depends smoothly on $\omega$ (resp. on $\omega$ and the origin of $\gamma$ on $\zeta$).
\end{remark}
The vectors $l$ and $n$ can be expressed as follows using the coordinate vector fields.
\begin{lemma}
\label{lncoord}
At each point of $D_T$, we have 
\begin{eqnarray*}
l=\frac{2}{N}\partial_v \, , ~ n=\frac{2}{N} \left( \partial_u+V_{\omega} \right) \, ,
\end{eqnarray*}
where $V_{\omega}$ is a vector field orthogonal to both $\nabla u$ and $\nabla v$ and s-smooth on $\mathfrak{D}_T$, $rV_\omega$ is smooth on $\mathfrak{D}_T$.
\end{lemma}

{\bf Proof.}
We first establish the result for $n$. Observing that $\nabla v$ is a null vector, the Ansatz $\nabla v=\alpha\partial_u+\beta\partial_v+W_{\omega}$, where $W_\omega$ is a vector field tangent to the $2$-surfaces of constant $u$ and $v$ (i.e. the surfaces orthogonal to both $\nabla u$ and $\nabla v$, in other words, the intersections of the forward and backward cones emanating from $\zeta$) gives~:
\begin{eqnarray*}
0=g(\nabla v, \nabla v)=\alpha \d v(\partial_u)+\beta \d v(\partial_v)+\gamma \d v (W_{\omega})=\beta.
\end{eqnarray*}
Now,
\[ \d u(\partial_u)=1\]
entails
\begin{eqnarray*}
\alpha=\d u( \nabla v)= \frac{2}{N} \d u (n) = \frac{2}{N} n^a \nabla_a u=\frac{4}{N^2}l_an^a=\frac{4}{N^2} \, .
\end{eqnarray*}
Hence
\[ n = \frac{N}{2} \nabla v = \frac{2}{N} \partial_u + \frac{N}{2} W_\omega \]
and we put $V_\omega := \frac{N^2}{4} W_\omega$. The vector fields $n$ and $\partial_u$ being s-smooth on $\mathfrak{D}_T$, so is $V_{\omega}$~; the same argument gives the regularity of $rV_\omega$ on $\mathfrak{D}_T$.

The case of $l$ is simpler. By construction, $\omega$ is constant on the integral lines of $\nabla u$ (and therefore of $l$) which are the null geodesics along the future null cones emanating from $\zeta$ (see lemma \ref{NablauGeod}). So is $u$ since $\nabla u$ is orthogonal and tangent to the future null cones from $\zeta$, i.e. the level hypersurfaces of $u$. Consequently, $l$ is colinear to $\partial_v$. Putting $l = \alpha \partial_v$,
\[ \d v \left( \partial_v \right) = 1 = \frac{1}{\alpha} \d v (l) = \frac{2}{N\alpha} g (n , l) = \frac{2}{N\alpha} \, .\]
This concludes the proof. \qed

\subsection{Structure at the tip of the cone viewed on $\mathfrak C$} \label{TipCone}

At $p_0$ (just as at any other point of the timelike curve $\zeta$), the Newman-Penrose tetrad which we have defined is multi-valued. We denote by $\gamma_{0,\omega}$ the integral curve of $\nabla u$ on ${\cal C}_0^+$ corresponding to the direction $\omega$ and by $(l_\omega , n_\omega , m_\omega , \bar{m}_\omega )$ the Newman-Penrose tetrad $(l,n,m,\bar{m})$ at $p_0$ corresponding to the direction $\omega$. Let $\{ o^A_\omega , \iota^A_\omega \}$ be the associated spin-frame. The vector $l_\omega$ points in a direction corresponding to $\omega$ and the vector $n_\omega$ points along another direction on $S^2$.
\begin{definition}
We denote by $\omega'$ the direction on $S^2$ corresponding to $n_\omega$ and call it the conjugate direction of $\omega$.
\end{definition}
As we shall see shortly, the pairs of conjugate directions play an important role in understanding the continuity at the vertex of a spinor field defined on the cone. The directions $\omega$ and $\omega'$ and their associated tetrads satisfy some important properties.
\begin{lemma} \label{RotationTetrad}
Given any direction $\omega \in S^2$ and $\omega '$ its conjugate direction, we have~:
\begin{itemize}
\item $(\omega')' = \omega$~;
\item the relation between the null tetrads $\{ l_\omega ,  n_\omega ,  m_\omega ,  \bar{m}_\omega \}$ and $\{ l_{\omega'} ,  n_{\omega '} ,  m_{\omega '} ,  \bar{m}_{\omega '} \}$ is given by
\[  l_{\omega'} = n_\omega \, ,~ n_{\omega '} = l_\omega  \, ,~ m_{\omega'} = e^{i\theta (\omega )} \bar{m}_{\omega } \, , ~ \bar{m}_{\omega '} = e^{-i \theta (\omega)} m_{\omega } \, ,\]
where $\theta (\omega ) \in \R / 2\pi \Z$ is a function of $\omega$ which is smooth on $S^2$ and such that
\[ \theta (\omega ) = \theta (\omega ' ) \, .\]
\end{itemize}
\end{lemma}
{\bf Proof.} The property that $(\omega')' = \omega$ follows from
\begin{equation} \label{AdaptedOmegaOmega'}
{\cal T}^a (p_0) = l^a_\omega + n^a_\omega = l^a_{\omega'} + n^a_{\omega'} \, .
\end{equation}
Indeed, any plane in the tangent space to $p_0$ contains at most two null directions. The plane spanned by $l^a_\omega$ and $n^a_\omega$ contains exactly two which are precisely $l^a_\omega$ and $n^a_\omega$. This plane also contains ${\cal T}^a (p_0 )$. Now by definition $l_{\omega'}$ is colinear to $n_\omega$ and $n_{\omega'} = {\cal T}^a (p_0 ) - l_{\omega'}$ is a null direction in this plane which is distinct from that of $l_{\omega'}$. Hence $n_{\omega'}$ must be colinear to $l_\omega$. Using \eqref{AdaptedOmegaOmega'} again, we get $l_{\omega'}^a = n_\omega^a $ and $n_{\omega '}^a = l_\omega^a$. This implies that the spin-frame is transformed as
\[ o_{\omega '}^A = \alpha \iota_{\omega }^A \, ,~  \iota_{\omega '}^A =  \beta o_{\omega}^A \]
where $\alpha, \beta \in \C$ and satisfy $\vert \alpha \vert = \vert \beta \vert =1$. Then,
\[ m_{\omega'}^a = o_{\omega '}^A  \bar{\iota}_{\omega '}^{A'} = \alpha \bar{\beta}  \iota_{\omega}^{A} \bar{o}_{\omega }^{A'}  = \alpha \bar{\beta} \bar{m}_\omega^a \, , \]
where $\vert \alpha \bar{\beta} \vert =1$. Putting $e^{i\theta} = \alpha \bar{\beta}$, the proof is complete. The other properties of $\theta$ follow by construction. \qed

We can establish a similar relation for the spin-frame modulo an overall sign choice.
\begin{lemma} \label{RotationDyad}
Given any direction $\omega \in S^2$ and $\omega '$ its conjugate direction, the relation between the dyads $\{ o_\omega^A ,  \iota_\omega^A \}$ and $\{ o_{\omega'}^A ,  \iota_{\omega '}^A \}$ is given by (after a choice of overall sign)
\[  o_{\omega'}^A = ie^{i\theta (\omega ) /2}\iota_\omega^A  \, ,~ \iota_{\omega'}^A = ie^{-i\theta (\omega )/2} o_{\omega }^A \, .\]
For the conjugate spin-frames, we have consequently
\[ \bar{o}_{\omega'}^{A'} = -ie^{-i\theta (\omega ) /2} \bar{\iota}_\omega^{A'}  \, ,~ \bar{\iota}_{\omega'}^{A'} = -ie^{i\theta (\omega )/2} \bar{o}_{\omega }^{A'} \, .\]

\end{lemma}
{\bf Proof.} We continue from the proof of Lemma \ref{RotationTetrad}. We must have
\[ 1 = \varepsilon_{AB} o_{\omega '}^A \iota_{\omega '}^B = \alpha \beta \varepsilon_{AB} \iota_{\omega}^A o_{\omega }^B = - \alpha \beta \]
and since $\vert \alpha \vert = \vert \beta \vert =1$, we must have
\[ \beta = -\frac{1}{\alpha} = -\bar{\alpha } \, .\]
Therefore,  $e^{i\theta} = \alpha \bar{\beta} = -\alpha^2$. Our choice of sign corresponds to $\alpha = ie^{i\theta /2}$. \qed

This entails the following relations between the components of a Dirac spinor at $p_0$ in the spin-frames $\{ o_\omega^A \, ,~ \iota_\omega^A \}$ and $\{ o_{\omega'}^A \, ,~ \iota_{\omega'}^A \}$~:
\begin{corollary} \label{RelationComponents}
Let $\Psi = \phi_A \oplus \chi^{A'}$ be a Dirac spinor at $p_0$, denote by $\Psi_i (\omega )$, $i=1,2,3,4$ its components (see \eqref{DiracSpinorComponents}) in the spin-frame $\{ o_\omega^A \, ,~ \iota_\omega^A \}$. We have
\begin{gather*}
\Psi_2 (\omega' ) = ie^{-i\theta (\omega ) /2} \Psi_1 (\omega) \, , ~\Psi_3 (\omega ') = ie^{i\theta (\omega ) /2} \Psi_4 (\omega) \, , \\
\Psi_1 (\omega') = ie^{i\theta (\omega )/2} \Psi_2 (\omega ) \, ,~ \Psi_4 (\omega ') = ie^{-i\theta (\omega )/2} \Psi_3 (\omega ) \, .
\end{gather*}
\end{corollary}
\begin{remark}
Applying this transformation twice leads to a global sign change of the components of $\Psi$. This is a typical consequence of the fact that the bundle of unitary spinor dyads is a two-fold covering of the bundle of normalized Newman-Penrose tetrads.
\end{remark}

\subsection{The geometry of the cone}

In this subsection, we derive a series of properties related to the geometry of the cone and define a Riemannian structure on the blown-up cone $\mathfrak{C}$.
\begin{proposition} \label{GaussCurvature}
The Gauss curvature $k$ of
\[ S_r := {\cal C}^+_0 \cap {\cal C}^-_{2r} \]
is equivalent to $1/r^2$ as $r \rightarrow 0$, meaning that $r^2 k \rightarrow 1$ as $r \rightarrow 0$.
\end{proposition}
{\bf Proof.} The surfaces $S_r$ are orthogonal to $l$ and $n$, consequently, the Gauss curvature of $S_r$ is given by (see \cite{PeRi}, Vol.1 p.272)
\[ k = 2 \Re \left( \sigma \sigma' - \Psi_2 - \rho \rho' + \Phi_{11} + \Lambda \right) \, ,\]
where $\Psi_2$ is the component of the Weyl spinor (encoding the information of the Weyl tensor which describes the conformal curvature) obtained by contracting it twice with $\iota^A$ and twice with $o^A$, $\Phi_{ab}$ is ($-1/2$ of) the trace-free part of the Ricci tensor,
\[ \Phi_{11} = \Phi_{ab}  l^a n^b \]
and 
\[ \Lambda = \frac{1}{24} \mathrm{Scal}_g \, ;\]
the spin-coefficients involved in the expression of the Gauss curvature are
\[ \sigma = m^a \nabla_m l_a \, ,~ \rho = m^a \nabla_{\bar{m}} l_a \]
and $\sigma'$ and $\rho'$ which are obtained by exchanging the roles of $o^A$ and $\iota^A$ in the expressions of $\sigma$ and $\rho$, i.e.
\[ \sigma' = \bar{m}^a \nabla_{\bar{m}} n_a \, ,~\rho' = \bar{m}^a \nabla_m n_a \, .\]
All the curvature components are bounded on $D_T$ since the metric is regular and the frame vectors and spinors are s-smooth.

Concerning $\rho$ and $\sigma$, if we work with the tetrad
\[ {\cal L} := \nabla u = \frac{2}{N} l \, ,~{\cal N} := \frac{N^2}{4} \nabla v = \frac{N}{2} n \, ,~ m \, ,~ \bar{m} \, ,\]
then we know from \cite{SeSchneEh} that (see appendix \ref{ConstraintEqSolve} for more details on the behaviour of $\rho$ and why this choice of tetrad is the right one for inferring such properties)
\[ \sigma \rightarrow 0 \mbox{ and } \rho \simeq -\frac{1}{r} \, .\]
Then a direct calculation shows that (recall that $N = \sqrt{2}$ at the tip of the cone)
\begin{eqnarray*}
m^a \nabla_m  l_a&=& \frac{N}{2} m^a \nabla_m {\cal L}_a \rightarrow 0 \mbox{ as } r \rightarrow 0 \, ,\\
m^a \nabla_{\bar{m}} l_a &=& \frac{N}{2} m^a \nabla_{\bar{m}} {\cal L}_a \simeq -\frac{1}{r\sqrt{2}} \mbox{ as } r \rightarrow 0 \, .
\end{eqnarray*}
Working with the tetrad
\[ \frac{N}{2} l \, ,~ \frac{2}{N} n = \nabla v \, ,~m \, ,~\bar{m} \, ,\]
we obtain the behaviour of $\sigma'$ and $\rho'$ near the tip of each past light-cone as follows (the change of sign for $\rho'$ is due to the fact that the cone is past-pointing and the vector $\nabla v$ future pointing)
\begin{eqnarray*}
\frac{2}{N} \sigma ' &=& \bar{m}^a \nabla_{\bar{m}}  \left( \frac{2}{N} n_a \right) \rightarrow 0 \mbox{ as } r \rightarrow 0 \, ,\\
\frac{2}{N} \rho ' &=& \bar{m}^a \nabla_{m} \left( \frac{2}{N} n_a \right) \simeq \frac{1}{r} \mbox{ as } r \rightarrow 0 \, .
\end{eqnarray*}
Hence
\[ \sigma' \rightarrow 0 \mbox{ and } \rho' \simeq \frac{1}{r\sqrt{2}} \mbox{ as } r \rightarrow 0 \]
It follows that
\[ k \simeq \frac{1}{r^2} \mbox{ as } r \rightarrow 0 \, . \qed\]
The Gauss-Bonnet formula implies that~:
\begin{corollary} \label{StrongRegularityCone}
The area of $S_r$ (its measure for the metric induced by $-g$, since the two induced metrics on $S_r$ are the same) behaves like $4\pi r^2$ near the vertex. Mapping $S_r$ onto the euclidian $2$-sphere via a diffeomorphism depending smoothly on $r\in [0,T]$, the measure induced by $-g$ on $S_r$, which is
\[ i m_a \bar{m}_b \, \d x^a \wedge \d x^b \]
is transformed into $r^2 \nu_r$ where $\nu_r$ is a measure on $S^2$ depending smoothly on $r\in [0,T]$ and uniformly equivalent to the euclidean measure on $S^2$. The $1$-forms $\frac{1}{r} m_a \d x^a$ and $\frac{1}{r} \bar{m}_a \d x^a$ are smooth (but not s-smooth) on $\mathfrak{C}$.
\end{corollary}
\begin{corollary} \label{ConeSmoothCylinder}
Let us consider on $\mathfrak{C}$ the metric
\[ \mathfrak{g}_{ab} \d x^a \d x^b := \d v^2 + \frac{2}{r^2} m_{(a} \bar{m}_{b)} \d x^a \d x^b \, . \]
This is a smooth Riemannian metric on $\mathfrak{C}$ which gives it the geometry of a smooth finite cylinder, part of its boundary being the $2$-sphere describing the blown up vertex of ${\cal C}^+_0$ (the other part is simply $S_{T}$). The vector fields
\[ l^a \, ,~ r m^a \, ,~ r\bar{m}^a \, ,\]
are smooth on $\mathfrak{C}$.
\end{corollary}
\begin{remark}
Note that the vector fields $r m^a$ and $r\bar{m}^a$ are just the duals for the metric $\mathfrak{g}$ of the $1$-forms $r^{-1} m_a$ and $r^{-1}\bar{m}_a$, i.e.
\[ r m^b \mathfrak{g}_{ab} = -\frac{1}{r} m_a \, ,~ r \bar{m}^b \mathfrak{g}_{ab} = -\frac{1}{r} \bar{m}_a \, .\]
\end{remark}

\section{Results} \label{Re}

\subsection{The $L^2$ setting}

Let $\Psi\in C^{\infty}(D_T;\S_A\oplus\S^{A'})$ be a smooth spinor field. The components $\Psi_1, \, \Psi_4$ are s-smooth functions on $\mathfrak{D}_T$. They therefore have s-smooth traces on $\mathfrak{C}$. We define ${\bf T}(\Psi)=(\Psi_1,\Psi_4)|_{\mathfrak{C}}$ to be these traces and 
\[ {\cal F}={\bf T}(C^{\infty}(D_T;\S_A\oplus\S^{A'})).\]

\begin{definition}
We define the space $L^2(({\cal C}_0^+, d\sigma_{{\cal C}^+_0});\C^2)$ as the completion of ${\cal F}$ in the norm
\[ \left\Vert (\Psi_1 , \Psi_4 ) \right\Vert^2_{L^2 ({\cal C}_0^+)} : = \int_{{\cal C}_0^+}(\vert\Psi_1\vert^2+\vert\Psi_4\vert^2) \d \sigma_{{\cal C}^+_0} \, ,\]
where $\d \sigma_{{\cal C}^+_0}=n \lrcorner \dvol^4$ and $\dvol^4$ is the $4$-volume measure induced by $g$. More explicitly,
\begin{gather*}
\dvol^4 = i l_a \d x^a \wedge n_b \d x^b \wedge m_c  \d x^c \wedge \bar{m}_d \d x^d \, , \\
\d \sigma_{{\cal C}^+_0} = i n_a \d x^a \wedge m_b \d x^b \wedge \bar{m}_c  \d x^c = i \frac{N}{2} \d v \wedge m_b \d x^b \wedge \bar{m}_c  \d x^c
\end{gather*}
and by proposition \ref{GaussCurvature} and corollary \ref{StrongRegularityCone}, the restriction of $\frac{i}{r^2} m_a \d x^a \wedge \bar{m}_b  \d x^b$ to the $2$-surface $S_r={\cal C}^+_0 \cap {\cal C}^-_{2r}$ is a smooth measure on the $2$-sphere which is bounded and bounded away from zero uniformly in $r \in ]0,T]$.
\end{definition}

Let $\Psi_T\in {\cal C}^\infty (\Sigma_T;\S_A\oplus \S^{A'})$. By the usual theorems 
for hyperbolic equations there exists a unique solution $\Psi = \phi_A\oplus\chi^{A'}
\in {\cal C}^\infty (D_T;\S_A\oplus\S^{A'})$ of (\ref{DirEq}) such that the trace of $\Psi$ on $\Sigma_T$ is equal to $\Psi_T$ (see \cite{Ni} for details). We can introduce the linear trace operator~:
\begin{definition}
Let $\Gamma$ be the operator which, to smooth data $\Psi_T \in {\cal C}^\infty (\Sigma_T;\S_A\oplus \S^{A'})$, associates the pair of complex scalar functions $(\Psi_1 , \Psi_4 )$ on $\mathfrak{C}$, first and fourth components in the spin-frame $(o^A , \iota^A )$ of the restriction of the corresponding solution $\Psi$ to the cone ${\cal C}^+_0$. By construction, we have
 \[ \Gamma:\, \Psi_T \in {\cal C}^\infty (\Sigma_T;\S_A\oplus \S^{A'}) \mapsto (\Psi_1,\Psi_4)\in {\cal F}\subset L^2(({\cal C}_0^+, d\sigma_{{\cal C}^+_0});\C^2). \]
\end{definition}
Using the conserved current we obtain by Stokes's theorem :
\[ \int_{\Sigma_T}*(\phi_A\bar{\phi}_{A'}\d x^{AA'}+\bar{\chi}_A\chi_{A'}\d x^{AA'})
=\int_{{\cal C}^+_0}*(\phi_A\bar{\phi}_{A'}\d x^{AA'}
+\bar{\chi}_A\chi_{A'}\d x^{AA'})\]
which can be written explicitly in terms of components of $\Psi$ as
\begin{equation}
\label{EqL2}
\int_{\Sigma_T} \left| \Psi \right|^2 \d \sigma_{\Sigma_T} = \int_{{\cal C}_0^+}(\vert\Psi_1\vert^2+\vert\Psi_4\vert^2) \d \sigma_{{\cal C}^+_0},
\end{equation}
where $\d \sigma_{\Sigma_T} = \frac{1}{2} {\cal T} \lrcorner \dvol^4$.

Equation \eqref{EqL2} entails that the operator $\Gamma$ possesses an extension to a bounded operator 
\begin{equation} \label{GammaL2}
\Gamma\in {\cal L} \left( L^2 \left( \Sigma_T;\S_A\oplus\S^{A'} \right) ;L^2 \left( ({\cal C}_0^+, d\sigma_{{\cal C}^+_0});\C^2 \right) \right) \, .
\end{equation}
Our first result is
\begin{theorem}
\label{T1}
The operator $\Gamma$ is an isometry.
\end{theorem}
The proof is given in section \ref{secProofT2}.

\begin{remark}
Any $2$-spinor $\phi_A$ at a point of ${\cal C}^+_0$ can be decomposed as
\[ \phi_A = \phi_1 o_A -\phi_0 \iota_A \, ,~ \phi_0 = \phi_A o^A \, ,~ \phi_1 = \phi_A \iota^A \, .\]
The spinor $o^A$ points along $l^a$ (in the sense of its flag-pole direction, see \cite{PeRi}, Vol. 1) which is tangent to ${\cal C}^+_0$, whereas $\iota^A$ points along $n^a$ which is transverse to ${\cal C}^+_0$. So we can consider $\phi_0$ as the part of $\phi_A$ transverse to ${\cal C}^+_0$ and $\phi_1$ as the part of $\phi_A$ tangent to ${\cal C}^+_0$. Similarly, for a Dirac spinor $\Psi$ at a point $p\in {\cal C}^+_0$, $(\Psi_1 , \Psi_4)$ is the part of $\Psi$ transverse to ${\cal C}^+_0$ and $(\Psi_2 , \Psi_3)$ the part tangent to ${\cal C}^+_0$. It is only the transverse part that appears on the r.h.s. of the energy equality \eqref{EqL2}.
\end{remark}

\subsection{Further $L^2$ estimates}

In this subsection we consider the Dirac equation with a source $\Xi = \rho_A \oplus \eta^{A'}$~:
\begin{equation} \label{DP}
\left\{ \begin{array}{ccl}{\left( \nabla^{AA'} -iq\Phi^{AA'} \right) \phi_A} &  = & {\frac{m}{\sqrt{2}} \chi^{A'}}+\eta^{A'}, \\
{\left( \nabla_{AA'} -iq\Phi_{AA'} \right)  \chi^{A'}} &  = & {-\frac{m}{\sqrt{2}} \phi_{A}}+\rho_A, \end{array} \right.
\end{equation}
We have~:
\begin{lemma}
\label{L1}
Let $\Xi = \rho_A\oplus\eta^{A'}\in {\cal C}^\infty (D_T;\S_A\oplus\S^{A'})$ and $\Psi = \phi_A\oplus\chi^{A'}\in {\cal C}^\infty (D_T;\S_A\oplus\S^{A'})$ be a smooth solution of \eqref{DP}. Then~:
\begin{equation} \label{L2E1}
\left| \int_{\Sigma_T} \left| \Psi \right|^2 \d \sigma_{\Sigma_T} - \int_{{\cal C}_0^+}(\vert\Psi_1\vert^2+\vert\Psi_4  \vert^2 ) d\sigma_{{\cal C}^+_0} \right| \lesssim \int_0^T\int_{\Sigma_t}( \vert \Psi \vert^2 + \vert \Xi \vert^2 ) \d \sigma_{\Sigma_t} \d t
\end{equation}
and for $0<t<T$,
\begin{equation} \label{L2E2}
\int_{\Sigma_t} \left| \Psi \right|^2 \d \sigma_{\Sigma_t} \lesssim \int_{\Sigma_T} \left| \Psi \right|^2 \d \sigma_{\Sigma_T} +\int_t^T\int_{\Sigma_s} \left| \Xi \right|^2 \d \sigma_{\Sigma_s} \d s \, .
\end{equation}
\end{lemma}
{\bf Proof.} We shall use the following notations for $0 \leq t_1 < t_2 \leq T$,
\begin{equation} \label{TimeIntervalCone}
{\cal S}_{t_1,t_2} := {\cal C}^+_0 \cap \{ t_1 \leq t \leq t_2 \} \, , ~ D_{t_1,t_2} := D_T \cap \{ t_1 \leq t \leq t_2 \} \, .
\end{equation}
It is sufficient to establish the two estimates for the Weyl equation~:
\begin{eqnarray}
\label{WP}
\nabla^{AA'}\phi_A=\eta^{A'}
\end{eqnarray}
with $\eta^{A'}\in {\cal C}^\infty (D_T;\S^{A'})$. We apply Stokes' theorem on the closed hypersurface made of $\Sigma_T$, $\Sigma_t$ and the $3$-surface ${\cal S}_{t,T}$. We obtain~:
\begin{eqnarray}
\label{EWL2}
\int_{\Sigma_T}*(\phi_A\bar{\phi}_{A'}dx^{AA'}) - \int_{\Sigma_t}*(\phi_A\bar{\phi}_{A'}dx^{AA'})=\int_{{\cal S}_{t,T}}\vert\phi_0\vert^2d\sigma_{{\cal C}^+_0}+2\Re\int_{D_{t,T}}\eta^{A'}\bar{\phi}_{A'}\dvol^4 \, .
\end{eqnarray}
We have~:
\begin{eqnarray*}
\left\vert\int_{D_{t,T}}\eta^{A'}\bar{\phi}_{A'}\dvol^4 \right\vert\lesssim\int_t^T\int_{\Sigma_s}\vert\eta^{A'}\bar{\phi}_{A'}\vert \d\sigma_{\Sigma_s}\d s\lesssim\int_t^T\int_{\Sigma_s}(T^{AA'}\phi_A\bar{\phi}_{A'}+T^{AA'}\bar{\eta}_A\eta_{A'})\d \sigma_{\Sigma_s}\d s \, .
\end{eqnarray*}
This entails
\[ \left| \int_{\Sigma_T} \left| \Psi \right|^2 \d \sigma_{\Sigma_T} - \int_{\Sigma_t} \left| \Psi \right|^2 \d \sigma_{\Sigma_t} - \int_{{\cal S}_{t,T}}(\vert\Psi_1\vert^2+\vert\Psi_4  \vert^2 ) d\sigma_{{\cal C}^+_0} \right| \lesssim \int_t^T\int_{\Sigma_s}( \vert \Psi \vert^2 + \vert \Xi \vert^2 ) \d \sigma_{\Sigma_s} \d s \]
which for $t=0$ gives (\ref{L2E1}) and (\ref{L2E2}) then follows via a Gronwall estimate. \qed

\subsection{Constraint equations on the cone}

The image of the data for \eqref{DirEq} by the trace operator $\Gamma$ only involves two components of the trace of the solution on ${\cal C}^+_{0}$~: $\Psi_1 = \phi_0$ and $\Psi_4 = -\chi_{0'}$. The other two are completely determined by the values of $\Psi_1$ and $\Psi_4$ on the cone and by the restriction of the Dirac equation tangent to ${\cal C}^+_{0}$~:
\begin{eqnarray}
\label{ConstraintEq}
 \left. \begin{array}{l}
{ l^\mathbf{a} (\partial_\mathbf{a}-iq\Phi_{\bf{a}})  \, \phi_1 - \bar{m}^\mathbf{a}
(\partial_\mathbf{a}-iq\Phi_{\bf{a}})  \, \phi_0 + (\alpha - \pi )\phi_0 + (\varepsilon -
\rho ) \phi_1 = - \frac{m}{\sqrt{2}} \chi_{0'} \, , } \\ \\
{ l^\mathbf{a} (\partial_\mathbf{a}-iq\Phi_{\bf{a}})  \, \chi_{1'} - m^\mathbf{a}
(\partial_\mathbf{a}-iq\Phi_{\bf{a}})  \, \chi_{0'} + (\bar{\alpha} - \bar{\pi} )\chi_{0'}
+ (\bar{\varepsilon} - \bar{\rho} ) \chi_{1'} = - \frac{m}{\sqrt{2}}
\phi_{0} \, . } \end{array} \right\}
\end{eqnarray}
On ${\cal C}^+_{0}$ we work with the coordinate system $(v,\omega )$. Using the expression $l = \frac{2}{N} \partial_v$ proved in lemma \ref{lncoord}, the equations \eqref{ConstraintEq} can be written as
\begin{eqnarray*}
\partial_v \phi_1 + \frac{N}{2} (\varepsilon - \rho ) \phi_1 &=& \frac{N}{2} \left( \bar{m}^a \partial_a - \alpha + \pi \right) \phi_0 + iq \frac{N}{2} \Phi_ a \left( l^a \phi_1 - \bar{m}^a \phi_0 \right) - \frac{Nm}{2\sqrt{2}} \chi_{0'} \, ,\\
\partial_v \chi_{1' }+ \frac{N}{2} (\bar\varepsilon - \bar\rho ) \chi_{1'} &=& \frac{N}{2} \left( {m}^a \partial_a - \bar\alpha + \bar\pi \right) \chi_{0'} + iq \frac{N}{2} \Phi_ a \left( l^a \chi_{1'} - {m}^a \chi_{0'} \right) - \frac{Nm}{2\sqrt{2}} \phi_{0} \, .
\end{eqnarray*}
It is a priori far from obvious that knowing $\Psi_1$ and $\Psi_4$ on the cone, we can solve these equations. A look at the form of the equations in the case of flat spacetime (given in appendix \ref{ConstraintEqSolve}) will convince the reader that there is a genuine difficulty there since they are singular at the vertex. The usual simplified understanding of the constraints for the Goursat problem is that the null data fix the value of the remaining components at the tip of the cone, in our case via the continuity matching conditions (see Corollary \ref{RelationComponents}),
\begin{equation} \label{MatchCond}
\left( \begin{array}{c} \Psi_2 \\ \Psi_3 \end{array} \right) (0,\omega ) = \left( \begin{array}{c} -ie^{-i\theta (\omega )/2} \Psi_1 \\ -ie^{i\theta (\omega )/2} \Psi_4 \end{array} \right) (0,\omega' ) \, ,
\end{equation}
and the restriction of the equation to the cone (our constraint equations \eqref{ConstraintEq}) allows to propagate these data from the tip and to recover the remaining part of the field on the cone. This is indeed true when one considers the Goursat problem on two intersecting null hypersurfaces. In the case of the cone, what really happens is the following~:
\begin{proposition} \label{WellPosedConstraintEq}
Given $(\Psi_1, \Psi_4 ) \in {\cal F}$, equations \eqref{ConstraintEq} have a unique solution $(\Psi_2, \Psi_3 )$ which is continuous on $\mathfrak{C}$ (the other solutions blow up at the tip) and it satisfies the matching conditions \eqref{MatchCond}.
\end{proposition}
The proof is long and technical~; it is given in appendix \ref{ConstraintEqSolve}.
\begin{definition}
We denote by $K$ the operator which to data $(\Psi_1 , \Psi_4 ) \in {\cal F}$ associates the solutions $(\Psi_2, \Psi_3 )$ of the constraint equations. It is a linear operator from $\cal F$ into $L^2 ({\cal C}^+_0)$.
\end{definition}
Moreover, as a straightforward consequence of the integration of the constraints detailed in appendix \ref{ConstraintEqSolve}, and more particularly of the expression of the rescaled constraint equations \eqref{RescConstraintEq}, we have~:
\begin{proposition} \label{RegularitySolConstraintEq}
For any $(\Psi_1 , \Psi_4 ) \in {\cal F}$, the corresponding solution $(\Psi_2 , \Psi_3 )$ of the constraints is of regularity ${\cal C}^1([0,2T]_v; C^{\infty}(S^2))$.
\end{proposition}

\subsection{The $H^1$ setting} \label{H1Setting}

Using the geodesic convexity of $\Omega$, we choose a global smooth orthonormal frame $\{ e_0 \, ,~ e_1 \, ,~ e_2 \, ,~ e_3 \}$ on $D_T$ such that $e_0 = \frac{\cal T}{ g ({\cal T} , {\cal T} )^{1/2}}$ and $e_\alpha$, $\alpha =1,2,3$ are consequently tangent to the hypersurfaces $\Sigma_t$. For $\alpha =0,1,2,3$, the spinor $\nabla_{e_\alpha} \Psi$ satisfies the equation (see \eqref{3+1Dirac2})~:
\begin{equation} \label{DirDalphaPsi}
( \notD + {\cal P} )(\nabla_{e_\alpha} \Psi ) = [ \notD +{\cal P}  , \nabla_{e_\alpha} ] \Psi
\end{equation}
and the error term $[ \notD +{\cal P}  , \nabla_{e_\alpha} ] \Psi$ is a smooth first order differential operator applied to $\Psi$ and is therefore controlled in norm for each $t$ by $\Vert \Psi \Vert_{L^2 (\Sigma_t )} + \sum_{\beta =0}^3 \Vert \nabla_{e_\beta} \Psi \Vert_{L^2 (\Sigma_t )}$. From lemma \ref{L1}, $\nabla_{e_\alpha} \Psi$ satisfies the estimates for $t \in [0,T]$~:
\begin{gather} 
\left| \Vert \nabla_{e_\alpha} \Psi \Vert^2_{L^2 (\Sigma_t)} - \int_{{\cal S}_{0,t}} ( \vert ( \nabla_{e_\alpha}\Psi )_1\vert^2+\vert (\nabla_{e_\alpha} \Psi)_4\vert^2) d\sigma_{{\cal C}^+_0} \right| \hspace{2in} \nonumber \\
\hspace{2in} \lesssim \int_0^t ( \Vert \Psi \Vert^2_{L^2 (\Sigma_s )} + \sum_{\beta =0}^3 \Vert \nabla_{e_\beta} \Psi \Vert^2_{L^2 (\Sigma_s )} ) \d s \, ,\label{EstBothWaysTPsi}
\end{gather}
\begin{equation} \label{EstTPsiH1}
\Vert \nabla_{e_\alpha} \Psi \Vert^2_{L^2 (\Sigma_t)} \lesssim \Vert \nabla_{e_\alpha} \Psi \Vert^2_{L^2 (\Sigma_T)} + \int_t^T ( \Vert \Psi \Vert^2_{L^2 (\Sigma_s )} + \sum_{\beta =0}^3 \Vert \nabla_{e_\beta} \Psi \Vert^2_{L^2 (\Sigma_s )} ) \d s \, .
\end{equation}
Using a Gronwall estimate, \eqref{EstBothWaysTPsi} and \eqref{EstTPsiH1} give the equivalence
\begin{equation} \label{H1NormEq}
\Vert \Psi \Vert^2_{L^2 (\Sigma_T )} + \sum_{\alpha =0}^3 \Vert \nabla_{e_\alpha} \Psi \Vert^2_{L^2 (\Sigma_T )} \simeq \int_{{\cal C}_0^+} ( \vert  \Psi_{1,4}\vert^2+\sum_{\alpha =0}^3 \vert (\nabla_{e_\alpha} \Psi )_{1,4} \vert^2 ) d\sigma_{{\cal C}^+_0} \, ,
\end{equation}
where $\vert  \Psi_{1,4}\vert^2$ denotes $\vert  \Psi_{1}\vert^2 + \vert  \Psi_{4}\vert^2$.

Using the equation (more particularly the form \eqref{3+1Dirac1}, \eqref{3+1Dirac2}), we see that the left-hand side of \eqref{H1NormEq} is equivalent to the $H^1$ norm of $\Psi$ on $\Sigma_T$. We now wish to understand the right-hand side as a function space involving purely $\Psi_1$ and $\Psi_4$ on ${\cal C}^+_0$. First we note that the right-hand side is equivalent to (this amounts to decomposing the complex vectors $(\nabla^a \phi_B ) o^B$ and $(\nabla^a \chi_{B'}) \bar{o}^{B'}$ on two different bases and using the equivalence of norms on $\C^4$)
\[ \int_{{\cal C}_0^+} ( \vert  \Psi_{1,4}\vert^2+ \vert (\nabla_{n} \Psi )_{1,4} \vert^2 + \vert (\nabla_{l} \Psi )_{1,4} \vert^2 + \vert (\nabla_{m} \Psi )_{1,4} \vert^2 + \vert (\nabla_{\bar{m}} \Psi )_{1,4} \vert^2 ) d\sigma_{{\cal C}^+_0} \, . \]
Hence the equivalence above becomes
\begin{equation} \label{H1NormEqTangent}
\Vert \Psi \Vert^2_{H^1 (\Sigma_T )} \simeq \int_{{\cal C}_0^+} ( \vert  \Psi_{1,4}\vert^2+ \vert (\nabla_{n} \Psi )_{1,4} \vert^2 + \vert (\nabla_{l} \Psi )_{1,4} \vert^2 + \vert (\nabla_{m} \Psi )_{1,4} \vert^2 + \vert (\nabla_{\bar{m}} \Psi )_{1,4} \vert^2 ) d\sigma_{{\cal C}^+_0} \, .
\end{equation}
The eight terms involving derivatives have the following expressions (given by a direct calculation using the Newman-Penrose formalism)~:
\begin{eqnarray*}
(\nabla_n \Psi )_1 &=& (\nabla_n \phi_A ) o^A = \nabla_n \phi_0 - \gamma \phi_0 + \tau \phi_1 = \nabla_n \Psi_1 - \gamma \Psi_1 + \tau \Psi_2 \, , \\
(\nabla_l \Psi )_1 &=& (\nabla_l \phi_A ) o^A = \nabla_l \phi_0 - \varepsilon \phi_0 + \kappa \phi_1 = \nabla_l \Psi_1 - \varepsilon \Psi_1 + \kappa \Psi_2 \, , \\
(\nabla_m \Psi )_1 &=& (\nabla_m \phi_A ) o^A = \nabla_m \phi_0 - \beta \phi_0 + \sigma \phi_1 = \nabla_m \Psi_1 - \beta \Psi_1 + \sigma \Psi_2 \, , \\
(\nabla_{\bar{m}} \Psi )_1 &=& (\nabla_{\bar{m}} \phi_A ) o^A = \nabla_{\bar{m}} \phi_0 - \alpha \phi_0 + \rho \phi_1 = \nabla_{\bar{m}} \Psi_1 - \alpha \Psi_1 + \rho \Psi_2 \, , \\
(\nabla_n \Psi )_4 &=& -(\nabla_n \chi_{A'} ) \bar{o}^{A'} = -\nabla_n \chi_{0'} + \bar\gamma \chi_{0'} - \bar\tau \chi_{1'} = \nabla_n \Psi_4 - \bar\gamma \Psi_4 - \bar\tau \Psi_3 \, , \\
(\nabla_l \Psi )_4 &=& -(\nabla_l \chi_{A'} ) \bar{o}^{A'} = -\nabla_l \chi_{0'} + \bar\varepsilon \chi_{0'} - \bar\kappa \chi_{1'} = \nabla_l \Psi_4 - \bar\varepsilon \Psi_4 - \bar\kappa \Psi_3 \, , \\
(\nabla_m \Psi )_4 &=& -(\nabla_m \chi_{A'} ) \bar{o}^{A'} = -\nabla_m \chi_{0'} + \bar\beta \chi_{0'} - \bar\sigma \chi_{1'} = \nabla_m \Psi_4 - \bar\beta \Psi_4 - \bar\sigma \Psi_3 \, , \\
(\nabla_{\bar{m}} \Psi )_4 &=& -(\nabla_{\bar{m}} \chi_{A'} ) \bar{o}^{A'} = -\nabla_{\bar{m}} \chi_{0'} + \bar\alpha \chi_{0'} - \bar\rho \chi_{1'} = \nabla_{\bar{m}} \Psi_4 - \bar\alpha \Psi_4 - \bar\rho \Psi_3 \, .
\end{eqnarray*}
The terms involving a transverse derivative (i.e. $(\nabla_n \Psi )_{1,4}$) can be expressed in terms of tangential derivatives~: using the first and third equations in \eqref{NewmanDirac} and remembering that on the cone
\[ \left( \begin{array}{c} {\Psi_2} \\ {\Psi_3} \end{array} \right) = K \left( \begin{array}{c} {\Psi_1} \\ {\Psi_4} 
\end{array} \right) \, ,\]
the expression becomes~:
\begin{eqnarray*}
\left( \begin{array}{c} {(\nabla_n \Psi)_1} \\ {(\nabla_n \Psi)_4} \end{array} \right)  &=& \left( \begin{array}{cc} { i q n^a \Phi_a - \mu} & 0 \\ 0 & {iqn^a \Phi_a - \bar\mu} \end{array} \right) \left( \begin{array}{c} {\Psi_1} \\ {\Psi_4} \end{array} \right) \\
&&+ \left( \begin{array}{cc} { \nabla_m -iqn^a \Phi_a + \beta} & {\frac{m}{\sqrt{2}}} \\ {-\frac{m}{\sqrt{2}}} & {-\nabla_{\bar{m}} +iq n^a \Phi_a - \bar\beta} \end{array} \right) K \left( \begin{array}{c} {\Psi_1} \\ {\Psi_4} 
\end{array} \right) \, .
\end{eqnarray*}
\begin{definition}
We denote by $L_n$ the operator on the right-hand side of the equation above, i.e.
\[ L _n= \left( \begin{array}{cc} { i q n^a \Phi_a - \mu} & 0 \\ 0 & {iqn^a \Phi_a - \bar\mu} \end{array} \right) + \left( \begin{array}{cc} { \nabla_m -iqm^a \Phi_a + \beta} & {\frac{m}{\sqrt{2}}} \\ {-\frac{m}{\sqrt{2}}} & {-\nabla_{\bar{m}} +iq \bar{m}^a \Phi_a - \bar\beta} \end{array} \right) K \, . \]
We also denote by $L_l$, $L_m$ and $L_{\bar{m}}$ the operators acting on $\Psi_{1,4}$ corresponding to the tangential derivatives (using the fact that $\kappa=0$)~:
\begin{eqnarray*}
L_l &:=& \nabla_l + \left( \begin{array}{cc} { -\varepsilon} & 0 \\ 0 & -\bar\varepsilon \end{array} \right)\, ,\\
L_m &:=& \nabla_m + \left( \begin{array}{cc} { -\beta} & 0 \\ 0 & {- \bar\beta} \end{array} \right) + \left( \begin{array}{cc} { \sigma} & 0 \\ 0 & -\bar\sigma \end{array} \right) K \, ,\\
L_{\bar{m}} &:=& \nabla_{\bar{m}} + \left( \begin{array}{cc} { -\alpha} & 0 \\ 0 & {- \bar\alpha} \end{array} \right) + \left( \begin{array}{cc} \rho & 0 \\ 0 & -{\bar\rho} \end{array} \right) K \, .
\end{eqnarray*}
\end{definition}
The main property which allows us to define our space of data associated to $H^1$ solutions is the following~:
\begin{proposition}
\label{Prop4.9}
The operators $L_n$, $L_l$, $L_m$ and $L_{\bar{m}}$ are well defined as operators from ${\cal F}$ to $L^2 ( ({\cal C}^+_0 \, ;~ \d \sigma_{C^+_0}) \, ;~\C^2 )$.
\end{proposition}
{\bf Proof.} The bulk of the work was done for the integration of the constraint equations on the cone in appendix \ref{ConstraintEqSolve} which lead to proposition \ref{RegularitySolConstraintEq}. Using the fact that $rm, r\bar{m}, l, n$ are smooth vector fields on $\mathfrak{C}$, all that remains to do now is to check that the spin-coefficients $\mu$, $\beta$, $\sigma$, $\alpha$, $\varepsilon$ and $\rho$ have a reasonable behaviour near the vertex. This is done in appendix \ref{MoreSpinCoeffs}. We find that once multiplied by $r$, they are all s-smooth on $\mathfrak{C}$ (see lemma \ref{SpinCoeffs}). This and proposition \ref{RegularitySolConstraintEq} give the result since the measure is $r^2$ times a measure on the $2$-sphere which varies smoothly with $r$ in $[0,T]$. \qed

\begin{remark}
Note that the spinor whose components are $(\Psi_1,...,\Psi_4)$ has $H^1$ regularity on the cone (see proof of Proposition \ref{WellPosedConstraintEq}). This entails the part of Proposition \ref{Prop4.9} concerning $L_l,\, L_m,\, L_{\bar{m}}$. $L_n$ however needs to be treated more carefully. The proof we gave of the proposition unifies the treatments of the four operators. 
\end{remark}

We now define on ${\cal C}^+_0$ the Hilbert space
\begin{definition} \label{H1DataSpace}
Let ${\cal H}_{{\cal C}^+_0}$ be the completion of ${\cal F}$ in the norm
\[ \left\| \left( \begin{array}{c} {\Psi_1} \\ {\Psi_4} \end{array} \right) \right\|^2_{ {\cal H}_{{\cal C}^+_0}} := \left\| \left( \begin{array}{c} {\Psi_1} \\ {\Psi_4} \end{array} \right) \right\|^2_{L^2 ( ({\mathfrak C} \, ;~ \d \sigma_{C^+_0}) \, ;~\C^2 )} + \sum_{\alpha \in \{ n,l,m,\bar{m} \}} \left\| L_\alpha  \left( \begin{array}{c} {\Psi_1} \\ {\Psi_4} \end{array} \right) \right\|^2_{L^2 ( ({\mathfrak C} \, ;~ \d \sigma_{C^+_0}) \, ;~\C^2 )} \, . \]
\end{definition}
Equivalence \eqref{H1NormEq} gives
\begin{lemma}
For all smooth data $\Psi_T \in {\cal C}^\infty (\Sigma_T;\S_A\oplus \S^{A'})$, we have
\begin{equation}
\label{H1est}
\Vert\Psi_T \Vert_{H^1(\Sigma_T)} \lesssim \Vert\Gamma\Psi_T \Vert_{{\cal H}_{{\cal C}^+_0}} \lesssim \Vert\Psi_T \Vert_{H^1(\Sigma_T)}
\end{equation}
and the trace operator $\Gamma$ therefore extends as a continuous operator from $H^1 (\Sigma_T )$ to ${\cal H}_{{\cal C}^+_0}$.
\end{lemma}

The main result of this paper, of which theorem \ref{T1} is a consequence, is~:
\begin{theorem} \label{T2}
$\Gamma$ is an isomorphism from $H^1 (\Sigma_T )$ onto ${\cal H}_{{\cal C}^+_0}$.
\end{theorem}
The proof is given in section \ref{secProofT2}.

\subsection{The Cauchy problem on a rough hypersurface on spatially compact spacetimes}

This section contains an extension to the Dirac equation in $4$ spacetime dimensions, of the results of \cite{Ho} for the Cauchy problem on a Lipschitz hypersurface for the wave equation. We consider a smooth compact manifold $X$ without boundary of dimension $3$. The spacetime ${\cal X} := \R_t \times X$ is endowed with a smooth Lorentzian metric $g$ of the form
\[ g= \frac{N^2}{2} \d t ^2 -h(t) \]
where $h(t)$ is a time-dependent Riemannian metric on $X$. We denote as before
\[ {\cal T} := \frac{2}{N} \frac{\partial}{\partial t} = N \nabla t \, .  \]
By \cite{Ge1, Ge2, Stie}, $\cal X$ admits a spin structure. We denote by $\vert \Psi \vert$ the norm induced by $\cal T$ on Dirac spinors at a point. Let $X_t$ denote the hypersurface $\{t \} \times X$ for any $t \in \R$. Using the parallelizability of $X$ (see \cite{Stie}), we consider on $\cal X$ a global smooth orthonormal frame
\[ e_0 = \frac{1}{\sqrt 2} {\cal T} \, ,~ e_\alpha \, ,~ \alpha = 1, 2, 3 \, .\]
Let $\cal S$ be a Cauchy hypersurface in $\cal X$ with low regularity, defined as the graph of a function $f \, : ~ X_0 \rightarrow \R$ which is merely assumed Lipschitz-continuous on $X$. We shall consider $f$ also as a function from $\cal X$ to $\R$ constant on the integral lines of $\partial_t$. Lipschitz continuous functions are differentiable almost everywhere, hence, the conormal to $\cal S$~:
\[ \nu_a \d x^a = \d t - \d f \]
is defined almost everywhere on $\cal S$ and is in $L^\infty ({\cal S})$. So are the normal vector field (normalized for convenience to have unit component along $e_0$)
\begin{equation}
V = \frac{N}{\sqrt{2}} \left( \nabla t - \nabla f \right) = e_0 - \frac{N}{\sqrt{2}} \nabla f = e_0 + \frac{N}{\sqrt{2}} \sum_{\alpha=1}^3 \nabla_{\alpha} f \,  e_\alpha  \, , ~ \nabla_\alpha f := e_\alpha  f = \d f (e_\alpha ) \, , \label{NormalS}
\end{equation}
and the tangent vector fields to $\cal S$
\[ \tau_\alpha = \frac{N}{\sqrt{2}} \nabla_{\alpha} f \, e_0 + e_\alpha \, . \]
Let us denote
\[ T_1 := \min_{x\in X} f(x) \, ,~T_2 := \max_{x\in X} f(x) \, .\]

\subsubsection{$L^2$ and $H^1$ estimates}

We assume that $\cal S$ is weakly spacelike, i.e. almost everywhere on $\cal X$, $g_{ab} V^a V^b \geq 0$, or equivalently,
\begin{equation} \label{WeakSpace}
\frac{N^2}{2} \sum_{\alpha =1}^3 \left( \nabla_{\alpha} f \right)^2 \leq 1 \, .
\end{equation}
We define the energy of a spinor field $\Psi = \phi_A \oplus \chi^{A'}$ as
\[ {\cal E}_{\cal S} (\Psi) := \int_{\cal S} * (\phi_A \bar\phi_{A'} + \bar\chi_A \chi_{A'} ) \d x^{AA'} \, .\]
This is always non-negative since $\cal S$ is weakly spacelike.

We obtain the following energy estimates for all smooth solutions $\Psi$ of the Dirac equation~:
\begin{proposition} \label{AprioriEstimatesS}
For any smooth solution $\Psi$ of \eqref{DirEq} on $\cal X$, for all $t\in \R$,
\[ \Vert \Psi \Vert^2_{L^2 (X_t )} = {\cal E}_{\cal S} (\Psi ) \, .\]
Moreover, there exists a constant $C>1$ depending only on $T_1$ and $T_2$ such that for any smooth solution $\Psi$ of \eqref{DirEq} on $\cal X$, for all $t \in [T_1 -1 , T_2+1]$,
\[ \frac{1}{C} \Vert \Psi \Vert^2_{H^1 (X_t ) } \leq \Vert \Psi \Vert^2_{H^1 (X_0 ) } \leq C \Vert \Psi \Vert^2_{H^1 (X_t ) } \]
and
\begin{equation} \label{EquivH1S}
\frac{1}{C} ( {\cal E}_{\cal S} (\Psi ) + {\cal E}_{\cal S} ( \nabla_{e_0} \Psi ) + \sum_{\alpha =1}^3 {\cal E}_{\cal S} ( \nabla_{\tau_\alpha} \Psi ) ) \leq \Vert \Psi \Vert^2_{H^1 (X_0 ) } \leq C ( {\cal E}_{\cal S} (\Psi ) + {\cal E}_{\cal S} ( \nabla_{e_0} \Psi ) + \sum_{\alpha =1}^3 {\cal E}_{\cal S} ( \nabla_{\tau_\alpha} \Psi ) ) \, .
\end{equation}
\end{proposition}
The proof is analogous to that of estimates \eqref{H1NormEq} and \eqref{H1NormEqTangent}.

\subsubsection{Cauchy problem}

We now assume that $\cal S$ is uniformly spacelike, i.e. there exists $0<\varepsilon <1$ such that, almost everywhere on $\cal X$, $g_{ab} V^a V^b \geq \varepsilon$, or equivalently,
\begin{equation} \label{UnifSpace}
\frac{N^2}{2} \sum_{\alpha =1}^3 \left( \nabla_{\alpha} f \right)^2 \leq 1-\varepsilon \, .
\end{equation}
In this case, we can get rid of the derivative along $e_0$ in estimate \eqref{EquivH1S}. To see this, we write the form \eqref{3+1Dirac1}, \eqref{3+1Dirac2} of the Dirac equation as follows
\begin{eqnarray*}
\left( e_0 - \sum_{\alpha=1}^3 \frac{N}{\sqrt{2}} \nabla_{\alpha} f \, e_\alpha \right) . \nabla_{e_0} \Psi &=& -\sum_{\alpha =1}^3 e_\alpha . \left( \nabla_{e_\alpha} + \frac{N}{\sqrt{2}} \nabla_{\alpha} f \, {\nabla}_{e_0} \right) \Psi  - {\cal P} \Psi \\
&=& - \sum_{\alpha =1}^3 e_\alpha . \nabla_{\tau_\alpha} \Psi - {\cal P} \Psi \, .
\end{eqnarray*}
Clifford multiplying by the vector
\[ W = e_0 - \sum_{\alpha=1}^3 \frac{N}{\sqrt{2}} \nabla_{\alpha} f \,  e_\alpha \, , \]
we obtain
\[ \left( 1 - \sum_{\alpha =1}^3 \frac{N^2}{2} \left| \nabla_{\alpha} f \, \right|^2 \right) \nabla_{e_0} \Psi = - W . \left[ \sum_{\alpha =1}^3 e_\alpha . \left( \nabla_{e_\alpha} +\frac{N}{\sqrt{2}} \nabla_{\alpha} f \, \nabla_{e_0} \right) \Psi  \right] - W . {\cal P} \Psi \, , \]
or equivalently
\[ \nabla_{e_0} \Psi = \frac{1}{1 - \sum_{\alpha =1}^3 \frac{N^2}{2} \left| \nabla_{\alpha} f \, \right|^2} W. \left( - \sum_{\alpha =1}^3 e_\alpha . \nabla_{\tau_\alpha}  \Psi  - {\cal P} \Psi \right) \, . \]
It follows
\[ \vert \nabla_{e_0} \Psi \vert^2 \leq \frac{C}{\varepsilon} ( \vert \psi \vert^2 + \sum_{\alpha =1}^3 \vert \nabla_{\tau_\alpha}  \Psi  \vert^2 )\, .\]
This and proposition \ref{AprioriEstimatesS} imply
\begin{lemma} \label{TraceOpS}
The trace operator
\[ \Gamma \, :~ {\cal C}^\infty (X_0 \, ;~\mathbb{S}_A \oplus \mathbb{S}^{A'} ) \longrightarrow L^2 ({\cal S} \, ;~\mathbb{S}_A \oplus \mathbb{S}^{A'}) \, ,\]
which to smooth data $\Phi$ on $X_0$ associates the trace on $\cal S$ of the smooth solution $\Psi$ of \eqref{DirEq} such that $\Psi \vert_{X_0} = \Phi$, extends as a continuous linear map still denoted $\Gamma$~:
\[ \Gamma \, :~ L^2 (X_0 \, ;~\mathbb{S}_A \oplus \mathbb{S}^{A'} ) \longrightarrow L^2 ({\cal S} \, ;~\mathbb{S}_A \oplus \mathbb{S}^{A'}) \, .\]
Moreover, $\Gamma$ satisfies for all $\Phi \in L^2 (X_0 \, ;~\mathbb{S}_A \oplus \mathbb{S}^{A'} )$,
\[ \Vert \Gamma \Phi \Vert^2_{L^2 ({\cal S} )} = \Vert \Phi \Vert^2_{L^2 (X_0 )} \, . \]
 This entails that $\Gamma$ is one-to-one and with closed range.

The restriction of $\Gamma$ to $H^1 (X_0 \, ;~\mathbb{S}_A \oplus \mathbb{S}^{A'} )$ is continuous from this space to $H^1 ({\cal S} \, ;~\mathbb{S}_A \oplus \mathbb{S}^{A'} )$, and satisfies
\[ \frac{1}{C} \Vert \Gamma \Phi \Vert^2_{H^1 ({\cal S} )} \leq \Vert \Phi \Vert^2_{H^1 (X_0 )} \leq \frac{C}{\varepsilon} \Vert \Gamma \Phi \Vert^2_{H^1 ({\cal S} )}\, . \]
\end{lemma}
\begin{remark}
Note that in the case of $H^1$ data, the solution is in $H^1_\mathrm{loc} ({\cal X})$ and $\Gamma$ is therefore a trace in the usual sense.
\end{remark}
We have the theorem~:
\begin{theorem} \label{CauchyLipschitz}
Let $\Phi \in L^2 ({\cal S})$, there exists a unique solution
\[ \Psi \in {\cal C} \left( \R_t \, ;~ L^2 (X) \right) \]
of \eqref{DirEq} such that
\[  \Psi \vert_{\cal S} = \Phi \, .\]
Moreover, if $\Phi \in H^1 ({\cal S})$, then
\[ \Psi \in {\cal C} \left( \R_t \, ;~ H^1 (X) \right) \cap {\cal C}^1 \left( \R_t \, ; ~L^2 (X) \right) \, . \]
\end{theorem}

\section{Proofs of the main results} \label{Pro}

\subsection{Proof of Theorem \ref{CauchyLipschitz}}

We consider a sequence of smooth hypersurfaces ${\cal S}_n$ approaching $\cal S$ as follows\footnote{For the existence of such a sequence of smooth hypersurfaces approaching $\cal S$, see \cite{Ho}, Lemma 3.}~: each ${\cal S}_n$ is defined as the graph of a smooth function $f_n \, :~X \rightarrow \R$, $f_n \rightarrow f$ in $L^\infty (X)$, $\d f_n \rightarrow \d f$ almost everywhere on $X$ and there exists $0<\delta <\varepsilon$ such that for each n
\begin{equation} \label{UnifSpaceSn}
\sum_{\alpha =1}^3 \frac{N^2}{2} \left( \nabla_{\alpha} f_n \right)^2 \leq 1- \delta \mbox{ almost everywhere on } {\cal S}_n \, ,
\end{equation}
which means in particular that the hypersurfaces ${\cal S}_n$ are spacelike uniformly in $x \in X$ and $n$.

Thanks to lemma \ref{TraceOpS}, all that we need to prove the theorem is to show that for data $\Phi \in H^1 ({\cal S})$, we can construct a solution $\Psi$ whose trace on $\cal S$ is $\Phi$, i.e. that the range of $\Gamma$ contains $H^1 ({\cal S})$ (which is dense in $L^2 ({\cal S})$). We push $\Phi$ along the flow of the vector field $\cal T$ as data $\Phi_n$ on ${\cal S}_n$. Since the sequence $\{ f_n \}_n$ is bounded in $W^{1,\infty} (X)$, not only is each $\Phi_n$ in $H^1 ({\cal S}_n \, ;~\mathbb{S}_A \oplus {\mathbb S}^{A'} )$, but the norm
\[ \Vert \Phi_n \Vert_{H^1 ({\cal S}_n \, ;~\mathbb{S}_A \oplus {\mathbb S}^{A'} )} \]
is bounded in $n$ (it would be constant if we defined the $H^1$ norm on each ${\cal S}_n$ as L. Hörmander did, as the $H^1$ norm on $X_0$ of the pull-back along the flow of $\cal T$, but our definition involves the metric at the points of the surfaces ${\cal S}_n$). By standard theorems, for each $n$, there exists a unique solution
\[ \Psi_n \in {\cal C} (\R_t \, ;~H^1(X)) \cap {\cal C}^1 (\R_t \, ;~L^2 (X)) \]
of \eqref{DirEq} such that $\Psi_n \vert_{{\cal S}_n} = \Phi_n$. Now by Lemma \ref{TraceOpS}, the sequence $\Psi_n$ is bounded in ${\cal C} (I \, ;~H^1(X)) \cap {\cal C}^1 (I \, ;~L^2 (X)) $ for any bounded time interval $I$ containing $0$ and such that $I\times X$ contains all hypersurfaces ${\cal S}_n$ and $\cal S$. Modulo the extraction of a subsequence, we can therefore assume that $\Psi_n$ converges weakly in $H^1 (I \times X)$ and in $H^1 (X_0)$, towards a solution
\[ \Psi \in {\cal C} (I \, ;~L^2 (X)) \, ,\]
of equation \eqref{DirEq} which naturally extends as
\[ \Psi \in {\cal C}  (\R_t \, ;~L^2 (X)) \, . \]
Since $\Psi (0) \in H^1 (X)$, it follows that $\Psi$ is more regular~:
\[ \Psi \in {\cal C} (\R_t \, ;~H^1(X)) \cap {\cal C}^1 (\R_t \, ;~L^2 (X)) \, . \]
Now, choosing $1/2 < s < 1$ and using the Rellich-Kondrachov compact embedding theorem, it follows that modulo the extraction of another subsequence, $\Psi_n$ converges towards $\Psi$ strongly in $H^{s} (I \times X )$, therefore by standard trace theorems, strongly in $L^2 ({\cal S})$.  It remains to prove that $\Gamma (\Psi (0)) = \Phi$, or more simply that the trace of $\Psi$ on $\cal S$ is equal to $\Phi$. To establish this last result, we project spinors on a given global spin-frame, still denoting $\Psi_n$, $\Psi$, $\Phi_n$ and $\Phi$ the vectors of the components of the correponding spinors in the spin-frame. We have
\begin{eqnarray*}
\int_X \vert \Psi (f(x),x) - \Phi_n (f_n(x),x ) \vert^2 \d \mu &=& \int_X \vert \Psi (f(x),x) - \Psi_n (f_n(x),x ) \vert^2 \d \mu \\
&\leq & 2\int_X \vert \Psi (f(x),x) - \Psi_n (f(x),x ) \vert^2 \d \mu \\
&&+ 2\int_X \vert \Psi_n (f(x),x ) - \Psi_n (f_n(x),x ) \vert^2 \d \mu \, ,
\end{eqnarray*}
where $\d \mu$ is the measure induced on $X$ by, say, $h(0)$.

The first integral on the right-hand side tends to zero since $\Psi_n \rightarrow \Psi$ strongly in $L^2 ({\cal S})$. As for the second, denoting $(f(x) , f_n (x))$ the interval between $f(x)$ and $f_n (x)$,
\begin{eqnarray*}
\int_X \vert \Psi_n (f_n(x),x ) - \Psi_n (f(x),x ) \vert^2 \d \mu &\leq &  \int_X \vert \int_{(f(x) , f_n (x))} \partial_t \Psi_n (t,x ) \d t \vert^2 \d \mu \\
&\leq & \int_X \vert f_n (x) - f(x) \vert \int_{(f(x) , f_n (x))} \vert \partial_t \Psi_n (t,x ) \vert^2 \d t \d \mu \\
&\leq & \sup_{x\in X} \vert f_n (x) -f(x) \vert \int_{I\times X} \vert \partial_t \Psi_n (t,x ) \vert^2 \d t \d \mu \, .
\end{eqnarray*}
The factor in front of the integral tends to zero since $f_n$ converges uniformly towards $f$ on $X$ and the integral is bounded since $\Psi_n$ is bounded in ${\cal C}^1 (I \, ;~L^2 (X))$. It follows that
\[ \int_X \vert \Psi (f(x),x) - \Phi_n (f_n(x),x ) \vert^2 \d \mu \]
tends to zero. But since by construction $\Phi_n (f_n(x),x )$ tends to $\Phi (f(x),x )$ uniformly on $X$, this implies that the trace of $\Psi$ on $\cal S$ is equal to $\Phi$. The proof is complete. \qed

\subsection{Proof of Theorems \ref{T1} and \ref{T2}}
\label{secProofT2}

For these proofs, we assume that our coordinate system and Newman-Penrose tetrad are defined on a subdomain of $\Omega$ that is slightly larger than $D$, namely on ${\cal J}^+ (\zeta (-\eta )) \cap {\cal J}^- (\zeta (2T+\eta ))$ for some $\eta>0$. This is always possible since $D$ is compact inside the open set $\Omega$.

First recall (see lemma \ref{lncoord}) that we have $l=\frac{2}{N} \partial_v,\, n=\frac{2}{N}(\partial_u+V_{\omega})$, where $V_{\omega}$ is an s-smooth vector field on our domain with $\zeta$ blown up as a cylinder. Also $V_{\omega}$ and $m$ lie in the tangent planes to the $2$-surfaces of constant $u$ and $v$, which means that $V^a_{\omega} \partial_a$ and $m^a\partial_a$ involve only derivatives with respect to $\omega$.

Using \eqref{NewmanDirac} we see that the Dirac equation takes the form~:
\begin{equation}
\label{DiracGamma}
\partial_t\Psi=i\tilde{H}\Psi \, ;  \tilde{H}=\gamma D_r+\tilde{P}_{\omega}+\tilde{Q} \, ,~ \gamma=Diag(1,-1,-1,1) \, ,
\end{equation}
where $D_r$ denotes $-i \partial_r$. Here $\tilde{P}_{\omega}$ is a differential operator with derivatives only in the angular directions and $\tilde{Q}$ is a potential. Note that the operators $\tilde{P}_{\omega}$ and $\tilde{Q}$ depend on $t$.

The equality \eqref{EqL2} and the inequality \eqref{H1est} show that the trace operators ($\Gamma$ considered as acting from $L^2 (\Sigma_T)$ to $L^2 ({\cal C}^+_0 \, ;~ \C^2 )$ or from $H^1 (\Sigma_T )$ to ${\cal H}_{{\cal C}_0^+}$) are injective with closed range. We therefore only have to show that ${\mathcal F}$ is contained in the range of the trace operators. We consider $g_1$ and $g_4$ in ${\cal F}$
and define $g_{2,3}$ by
\begin{equation} \label{g23}
\left(\begin{array}{c} g_2\vert_{{\cal C}^+_0} \\ g_3\vert_{{\cal C}^+_0} \end{array}\right) := K \left( \begin{array}{c} g_1\vert_{{\cal C}^+_0} \\  g_4\vert_{{\cal C}^+_0} \end{array} \right) \, .
\end{equation}
By proposition \ref{RegularitySolConstraintEq}, $g_{2,3}\in C^1([0,2T];C^{\infty}(S^2))$.

Let us now open the cone by a factor $0<\lambda<1,\, \vert\lambda-1\vert<<1$. The new cone is
\[ {\cal C}_0^{+,\lambda}=\{(\lambda r,r,\omega);\, 0\le r\le T/\lambda,\, \omega\in S^2\} \, . \]
The tangent plane to ${\cal C}_0^{+,\lambda}$ at a given point $p$ is given for $r\neq 0$ by
\begin{eqnarray*}
T_p {\cal C}_0^{+,\lambda}&=&Span\{\lambda\partial_t+\partial_r,m, \bar{m} \}=Span\{(1+\lambda)l+(\lambda-1)n+(1-\lambda) V_{\omega},m, \bar{m}\}\\
&=&Span\{\tau_\lambda := (1+\lambda)l+(\lambda-1)n,m,\bar{m}\}.
\end{eqnarray*}

The cone ${\cal C}_0^{+,\lambda}$ is spacelike and we can therefore solve the corresponding Cauchy problem by theorem \ref{CauchyLipschitz}. We recover the solution of the Goursat problem in the limit $\lambda\rightarrow 1$.  
More precisely, for $\lambda<1$, we extend $g_{1,4}$ to s-smooth functions on the cone ${\cal C}_0^{+,\lambda}$ blown up
and $g_{2,3}$ to solutions of \eqref{g23} up to $v = 2 T/\lambda$~; $g_{2,3}$ will thus belong to $C^1([0,2 T/\lambda];C^{\infty}(S^2))$. We consider the Cauchy problem~:
\begin{eqnarray}
\label{CPClambda}
\left\{\begin{array}{rcl} \partial_t\Psi^{\lambda}&=&i\tilde{H}\Psi^{\lambda},\\
\Psi^{\lambda}(\lambda r,r,\omega)&=&g(2r,\omega);\, (r,\omega)\in [0,\frac{T}{\lambda}]\times S^2. \end{array}\right.
\end{eqnarray}
The quadruplet of functions $g$ are the components of a Dirac spinor on ${\cal C}_0^{+,\lambda}$ and thanks to the proof of Proposition \ref{WellPosedConstraintEq} this Dirac spinor is in $H^1({\cal C}_0^{+,\lambda})$.
By theorem \ref{CauchyLipschitz} the problem \eqref{CPClambda} has a unique solution with values in $H^1$ on the slices and it satisfies the estimates
\begin{gather}
\label{40a}
\Vert \Psi^\lambda \Vert^2_{L^2 (\Sigma_{T}^\lambda )} = {\cal E}_{{\cal C}^{+,\lambda}_0} (\Psi^\lambda ) \, , \\
\Vert \Psi^\lambda \Vert^2_{H^1 (\Sigma_T^\lambda ) } \lesssim {\cal E}_{{\cal C}^{+,\lambda}_0} (\Psi^\lambda ) + {\cal E}_{{\cal C}^{+,\lambda}_0} ( \nabla_{\tau_\lambda} \Psi^\lambda ) + {\cal E}_{{\cal C}^{+,\lambda}_0} ( \nabla_{m} \Psi^\lambda ) + {\cal E}_{{\cal C}^{+,\lambda}_0} ( \nabla_{\bar{m}} \Psi^\lambda )  + {\cal E}_{{\cal C}^{+,\lambda}_0} ( \nabla_{e_0} \Psi^\lambda ) \, .
\end{gather}
By construction the first four terms on the right-hand side are bounded uniformly in $\lambda$. It remains to estimate the last term which is equivalent (with constants uniform in $\lambda$ and $\Psi$) to ${\cal E}_{{\cal C}^{+,\lambda}_0} ( \nabla_{\cal T} \Psi^\lambda )$. For this, we give a precise expression of the energy of a spinor on ${\cal C}^{+,\lambda}_0$. Note that 
\[l_{\lambda}=\frac{1}{2}(1+\frac{1}{\lambda})l+\frac{1}{2}(\frac{1}{\lambda}-1)n\]
is orthogonal to ${\cal C}_0^{+,\lambda}$ and
\[n_{\lambda}=\frac{1}{2}(1+\lambda)n+\frac{1}{2}(\lambda-1)l\]
is transverse to ${\cal C}_0^{+,\lambda}$. We also have~:
\[g(l_{\lambda},n_{\lambda})=1.\]
Therefore we obtain~:
\[ {\cal E}_{{\cal C}^{+,\lambda}_0} (\Phi )=\int_{{\cal C}_0^{+,\lambda}}\left(\frac{1+\lambda}{2\lambda}\vert\Phi_{1,4}\vert^2+\frac{1-\lambda}{2\lambda}\vert\Phi_{2,3}\vert^2\right) \d \sigma_{{\cal C}_0^{+,\lambda}}\]
with $\d \sigma_{{\cal C}_0^{+,\lambda}}=n_{\lambda}\lrcorner \dvol^4$.

We now express the four components of $ \nabla_{\cal T} \Psi^\lambda $. Since
\[ \nabla_{\cal T} = l+n = \frac{2}{N} ( \partial_u + \partial_v + V_\omega ) = \frac{2}{N} ( \partial_t + V_\omega ) \]
for the components $1$ and $4$ we can use the calculations done in section \ref{H1Setting}~; for the two other components, we perform similar calculations~:
\begin{eqnarray*}
( \nabla_{\cal T} \Psi^\lambda )_1 &=& \frac{2}{N} \frac{\partial \Psi^\lambda_1}{\partial t} - ( \gamma + \varepsilon ) \Psi^\lambda_1 + ( \kappa + \tau ) \Psi^\lambda_2 + \frac{2}{N} V_\omega \Psi_1^\lambda\, ,\\
( \nabla_{\cal T} \Psi^\lambda )_2 &=& \frac{2}{N} \frac{\partial \Psi^\lambda_2}{\partial t} - ( \pi + \nu ) \Psi^\lambda_1 + ( \gamma + \varepsilon ) \Psi^\lambda_2 + \frac{2}{N} V_\omega \Psi_2^\lambda\, , \\
( \nabla_{\cal T} \Psi^\lambda )_3 &=& \frac{2}{N} \frac{\partial \Psi^\lambda_3}{\partial t} + ( \bar\pi + \bar\nu ) \Psi^\lambda_4 + ( \bar\gamma + \bar\varepsilon ) \Psi^\lambda_3 + \frac{2}{N} V_\omega \Psi_3^\lambda\, , \\
( \nabla_{\cal T} \Psi^\lambda )_4 &=& \frac{2}{N} \frac{\partial \Psi^\lambda_4}{\partial t} - ( \bar\gamma + \bar\varepsilon ) \Psi^\lambda_4 - ( \bar\kappa + \bar\tau ) \Psi^\lambda_3 + \frac{2}{N} V_\omega \Psi_4^\lambda\, .
\end{eqnarray*}
The vector field $V_\omega$ is s-smooth on $D_T$ and is in the span of $m$ and $\bar{m}$, so just as for $\tau_\lambda$, $m$ and $\bar{m}$, the  energy on ${\cal C}^+_0$ of $\nabla_{V_\omega} \Psi^\lambda$ is controlled uniformly in $\lambda$. Hence, putting $\Phi^{\lambda}=\partial_t \Psi^\lambda$ and using the behaviour of the spin-coefficients obtained in appendix \ref{MoreSpinCoeffs}, it is sufficient to estimate
\[ \int_{{\cal C}_0^{+,\lambda}}\left(\frac{1+\lambda}{2\lambda}\vert\Phi^\lambda_{1,4}\vert^2+\frac{1-\lambda}{2\lambda}\vert\Phi^\lambda_{2,3}\vert^2\right) \d \sigma_{{\cal C}_0^{+,\lambda}} \, . \]
Therefore we have to calculate $\Phi^{\lambda}(\lambda r,r,\omega)$. To this purpose we introduce the following coordinates~:
\begin{eqnarray*}
\left. \begin{array}{rcl} \tau&=&t-\lambda r \, , \\ x&=&r \end{array}
\right\}\Rightarrow \partial_t=\partial_{\tau} \, ;\, \partial_{r}=\partial_x-\lambda\partial_{\tau}\, .
\end{eqnarray*}
We have 
\begin{eqnarray}
\label{VII.18}
\partial_t\Psi^{\lambda}&=&i\tilde{H}\Psi^{\lambda}\nonumber\\
\Leftrightarrow \partial_{\tau}\Psi^{\lambda}&=&\left(1+\gamma\lambda\right)^{-1}
\left(\gamma\partial_x\Psi^{\lambda}+i(\tilde{P}_{\omega}+\tilde{Q})\Psi^{\lambda}\right).
\end{eqnarray}
Using (\ref{VII.18}) we calculate~:
\begin{eqnarray*}
\frac{1}{i} \partial_t \Psi^{\lambda}&=&Diag\left(\frac{1}{1+\lambda}, \frac{1}{\lambda-1}, \frac{1}{\lambda-1},\frac{1}{1+\lambda}\right)D_x\Psi^{\lambda}\\
&+&Diag\left(\frac{1}{1+\lambda}, \frac{1}{1-\lambda}, \frac{1}{1-\lambda},\frac{1}{1+\lambda}\right)\left(\tilde{P}_{\omega}+\tilde{Q}\right)\Psi^{\lambda} \, .
\end{eqnarray*}
Recalling that $g(2r,\omega)=\Psi^{\lambda}(\lambda r, r,\omega)$
we find :
\begin{eqnarray}
\label{VII.19}
\frac{1}{i} \Phi^{\lambda}(\lambda r, r,\omega)&=& 2 Diag\left(\frac{1}{1+\lambda}, \frac{1}{\lambda-1}, \frac{1}{\lambda-1},\frac{1}{1+\lambda}\right) D_{v}g (2r,\omega)\nonumber\\
&+&Diag\left(\frac{1}{1+\lambda}, \frac{1}{1-\lambda}, \frac{1}{1-\lambda},\frac{1}{1+\lambda}\right) \left(\tilde{P}_{\omega}+\tilde{Q}\right)g (2r,\omega)\nonumber\\
&=:&g^{\lambda}_H \, .
\end{eqnarray}
Note that
\begin{equation} \label{g23=0}
(g^{\lambda}_H)_{2,3}=0,
\end{equation}
because $g$ satisfies the constraint equations along the cone~:
\[ 2\partial_vg_{2,3}=i((\tilde{P}_{\omega}+\tilde{Q})g)_{2,3}.\]
It follows that ${\cal E}_{{\cal C}^{+,\lambda}_0} ( \nabla_{e_0} \Psi^\lambda )$ is uniformly bounded in $\lambda$ for $\vert\lambda-1\vert<<1.$ Using \eqref{40a} we see that we have uniformly in $\lambda$~:
\[\Vert\Psi^{\lambda}(T,.)\Vert_{H^1(\Sigma_T)}\lesssim 1. \]
Repeating the above arguments for the spaces $H^1(\Sigma_t)$
we see that we can extract a subsequence, still denoted $\Psi^{\lambda}$, s.t.
\begin{eqnarray*}
\Psi^{\lambda}\rightharpoonup \Psi &\quad &H^1(\Sigma_T) \, , \\
\Psi^{\lambda}\rightharpoonup\Psi&\quad & H^1(D_T) \, ,\\
\Psi^{\lambda}\rightarrow\Psi &\quad & H^{s}(D_T)\quad \forall\, 1/2<s<1. \,
\end{eqnarray*}
$\Psi$ is a solution of the Dirac equation and we have :
\begin{eqnarray}
\label{VII.16b}
||\Psi||_{H^1(\Sigma_T)}\lesssim ||g||_{{\cal H}_{{\cal C}^+_0}},\, ||\Psi||_{H^1(D_T)}\lesssim ||g||_{{\cal H}_{{\cal C}^+_0}}.
\end{eqnarray}
We want to check that 
\begin{eqnarray*}
\Psi_{1,4}(r,r,\omega)&=&g_{1,4}(2r,\omega)\quad \forall 0\le r\le T.
\end{eqnarray*}
In fact we can even show :
\[ \Psi(r,r,\omega)=g(2r,\omega). \]
We estimate :
\begin{eqnarray*}
\lefteqn{\int_{0}^{T}\int_{S^2}|g(2r,\omega)-\Psi^{\lambda}(r,r,\omega)|^2 \d \sigma_{{\cal C}_0^+}}\\
&=&\int_{0}^{T}\int_{S^2}|\Psi^{\lambda}(\lambda r,r,\omega)-\Psi^{\lambda}(r,r,\omega)|^2 \d \sigma_{{\cal C}_0^+}\\
&=&\int_{0}^{T}\int_{S^2}\left|\int_{\lambda r}^{r}\partial_t\Psi^{\lambda}(t,r,\omega)\d t\right|^2 \d \sigma_{{\cal C}_0^+}\\
&\le&|\lambda-1|T\int_{0}^{T}\int_{S^2}\int_{\lambda r}^{r}|\partial_t\Psi^{\lambda}(t,r,\omega)|^2 \d t \, \d \sigma_{{\cal C}_0^+} \\
&\lesssim&T^2|\lambda-1| \left( \Vert \tilde{H}\Psi^{\lambda}\Vert_{L^2({\Sigma}_T^\lambda)}^2 + \Vert\Psi^{\lambda}\Vert_{L^2({\Sigma}_T^\lambda)}^2 \right)
\lesssim T^2|\lambda-1|\rightarrow 0.
\end{eqnarray*}
On the other hand :
\begin{eqnarray*}
\int_{0}^{T}\int_{S^2}|\Psi^{\lambda}(r,r,\omega)-\Psi(r,r,\omega)|^2 \d \sigma_{{\cal C}_0^+}
\le||\Psi^{\lambda}-\Psi||^2_{H^{s}(D_T)}\rightarrow 0.
\end{eqnarray*}
Thus $\Psi(r,r,\omega)=g(2r,\omega)$. \qed

\appendix

\section{Solving the constraint equations on the cone} \label{ConstraintEqSolve}

\subsection{The flat case}

We work on Minkowski's spacetime $\mathbb{M}$ with spherical coordinates $(t,r,\theta ,\varphi )$. The Minkowski metric is given by
\[ \eta = \d t^2 - \d r^2 - r^2 \d \omega^2 \, ,~ \d \omega^2 = \d \theta^2 + \sin^2 \theta \, \d \varphi^2 \, .\]
The reference timelike curve used to construct the double null foliation will be the $\{ r=0 \}$ line. This gives us the standard null coordinates
\[ u = t-r \, ,~ v = t+r \, .\]
The natural associated null tetrad is
\[ l^a \partial_a = \frac{1}{\sqrt{2}} \left( \partial_t + \partial_r \right) \, ,~ n^a \partial_a = \frac{1}{\sqrt{2}} \left( \partial_t - \partial_r \right) \, ,~ m^a \partial_a = \frac{1}{r\sqrt{2}} \left( \partial_\theta + \frac{i}{\sin \theta} \partial_\varphi \right) \, .\]
The decomposition of the Dirac equation (charged or not is irrelevent here since the spacetime has no charge for the field to interract with) in the corresponding spin-frame is as follows (see for example \cite{Ch})
\begin{eqnarray*}
\left( \partial_t - \partial_r - \frac{1}{r} \right) \phi_0 - \frac{1}{r} \left( \partial_\theta + \frac12 \cot \theta + \frac{i}{\sin \theta} \partial_\varphi \right) \phi_1 &=& m \chi_{1'} \, , \\
\left( \partial_t + \partial_r + \frac{1}{r} \right) \phi_1 - \frac{1}{r} \left( \partial_\theta + \frac12 \cot \theta - \frac{i}{\sin \theta} \partial_\varphi \right) \phi_0 &=& -m \chi_{0'} \, , \\
\left( \partial_t - \partial_r - \frac{1}{r} \right) \chi_{0'} - \frac{1}{r} \left( \partial_\theta + \frac12 \cot \theta - \frac{i}{\sin \theta} \partial_\varphi \right) \chi_{1'} &=& m \phi_{1} \, , \\
\left( \partial_t + \partial_r + \frac{1}{r} \right) \chi_{1' }- \frac{1}{r} \left( \partial_\theta + \frac12 \cot \theta + \frac{i}{\sin \theta} \partial_\varphi \right) \chi_{0'} &=& -m \phi_{0} \, .
\end{eqnarray*}
The constraint equations along the forward cone $\{ t=r \}$ are the part of the Dirac equation which is tangent to the cone, i.e. the second and fourth equations above. We shall ignore the singularities involving $\sin \theta$, they are dealt with in the usual way by picking two well-chosen charts on the sphere and working in each chart with spherical coordinates for which the North and South poles are not contained in the chart. The presence of $1/r$ in several parts of the equation is a more fundamental difficulty. In the very simple situation we are considering, there is an explicit way of getting rid of these terms. We simply multiply the equations by $r$ and put
\[ \hat\phi_0 := \phi_0 \, ,~ \hat\phi_1 := r \phi_1 \, ,~ \hat\chi_{0'} := \chi_{0'} \, ,~ \hat\chi_{1'} := r \chi_{1'} \, .\]
The constraint equations along the cone then become
\begin{eqnarray*}
\left( \partial_t + \partial_r  \right) \hat\phi_1 - \left( \partial_\theta + \frac12 \cot \theta - \frac{i}{\sin \theta} \partial_\varphi \right) \hat\phi_0 &=& -mr \hat\chi_{0'} \, , \\
\left( \partial_t + \partial_r  \right) \hat\chi_{1' }- \left( \partial_\theta + \frac12 \cot \theta + \frac{i}{\sin \theta} \partial_\varphi \right) \hat\chi_{0'} &=& -mr \hat\phi_{0} \, .
\end{eqnarray*}
which are now clearly (modulo the remarks above) integrable from the tip of the cone. We can choose any initial data we wish for $\hat\phi_1$ and $\hat\chi_{1'}$ and this will give rise to a solution $(\phi_1 , \chi_{1'})$ of the constraint equations along the cone, but only the one corresponding to zero initial data for $\hat\phi_1$ and $\hat\chi_{1'}$ will be bounded at the vertex. The boundedness and even the continuity at the blown-up vertex of the corresponding solution, and the exact way in which its values there are determined by those of $\phi_0$ and $\chi_{0'}$, can be seen in the general case using the Newman-Penrose formalism (or rather its compacted version, also referred to as the Geroch-Held-Penrose formalism). The proof is not made easier by considering the case of flat spacetime, so we stop here the example of flat spacetime and move on to the proof in the general case.

\subsection{General case}

The method is the same as in the flat case but is described differently. The multiplication by $r$ of two components out of four will be here understood as scaling by a factor $r$ the spinor $\iota^A$ and keeping $o^A$ unchanged. The behaviour of spinor components and of the different parts of the equation under rescalings of the frame spinors (including, as is the case here, rescalings which do not preserve unitarity) is best described using the compacted spin-coefficient formalism (see Penrose and Rindler \cite{PeRi} Vol 1, section 4.12). We give a quick description of this formalism in appendix \ref{CompactedNP}.

The constraint equations along the cone are
\begin{equation} \label{ConstraintCSCF}
 \left. \begin{array}{l}
 { \mbox{\th} \phi_1 - \eth' \phi_0 - \pi \phi_0 - \rho \phi_1 -iq l^a \Phi_a \phi_1 + i q \bar{m}^a \Phi_a \phi_0 = - \frac{m}{\sqrt{2}} \chi_{0'} \, , } \\ \\
{ \mbox{\th} \chi_{1'} - \eth \chi_{0'} - \bar{\pi} \chi_{0'} - \bar{\rho} \chi_{1'} -iq l^a \Phi_a \chi_{1'} + i q m^a \Phi_a \chi_{0'} = - \frac{m}{\sqrt{2}}
\phi_{0} \, .} \end{array} \right\}
\end{equation}
The behaviour of the solution to the constraint equation at the vertex will be determined by the behaviour of some spin-coefficients there, in particular of the geodesic convergence $\rho$, which is known provided the vector $l$ of our null tetrad on the cone is a gradient field. So we work on ${\cal C}^+_0$ with the null tetrad~:
\begin{equation}
{\cal L} := \nabla u = \frac{2}{N} l \, ,~ {\cal N} := \frac{N^2}{4} \nabla v = \frac{N}{2} n \, ,~ m \, ,~ \bar{m}
\end{equation}
and consider the equations \eqref{ConstraintCSCF} in reference to this null tetrad and the associated spin-frame
\[ \{ O^A \, ,~ I^A \} = \left\{ \sqrt{\frac{2}{N}} o^A \, ,~ \sqrt{\frac{N}{2}} \iota^A \right\} \]
instead of $\{ l^a , n^a , m^a , \bar{m}^a \}$ and $\{ o^A \, ,~ \iota^A \}$.
\begin{remark}
Note that on the cone we have
\[ {\cal L} = \frac{4}{N^2} \frac{\partial}{\partial v} \]
as a direct consequence of lemma \ref{lncoord}.
\end{remark}

This choice of null tetrad guarantees that $\cal L$ is a gradient field but also (the following result is established in the more complete proposition (7.1.60) in Penrose and Rindler Vol 2 \cite{PeRi}, we give a direct proof here)
\begin{lemma} \label{NablauGeod}
The integral curves of $\cal L$ are geodesics and $\cal L$ corresponds to a choice of affine parameter, i.e.
\[\nabla_{\cal L} {\cal L} = 0 \, .\]
\end{lemma}
{\bf Proof.} The proof is direct~:
\begin{eqnarray*}
\nabla_{\cal L} {\cal L}^b &=& \nabla_{\nabla u} {\nabla^b u} \, , \\
&=& \nabla_a u \nabla^a \nabla^b u \, ,\\
&=& \nabla_a u \nabla^b \nabla^a u \mbox{ since the connection is torsion-free,} \\
&=& \nabla^b \left( \nabla_a u \nabla^a u \right) - \left( \nabla^b \nabla_a u \right) \nabla^a u \, ,\\
&=& 0 - \nabla_a u \nabla^a \nabla^b u \mbox{ since } \nabla u \mbox{ is null and the connection torsion-free,} \\
&=& - \nabla_{\nabla u} \nabla^b u  \, . \qed
\end{eqnarray*}
An important coefficient describing the geodesic congruence generated by $\cal L$ is the spin-coefficient
\[ \rho = m^a \nabla_{\bar{m}} {\cal L}_a \, .\]
It is referred to as the geodesic convergence of the congruence which is justified by the following result.
\begin{lemma}
The commutator of $m$ and $\bar{m}$ is orthogonal to both $\cal L$ and $\cal N$. This together with $\nabla_{\cal L} {\cal L} =0$ implies in particular that
\[ \rho = -\frac{1}{2} \mathrm{div} {\cal L} \, .\]
\end{lemma}
\begin{remark}
Note that this is established in Penrose and Rindler \cite{PeRi}, it follows from Vol. 1 equation (5.12.13) and Vol. 2 equation (7.1.27). We give a proof here for the convenience of the reader.
\end{remark}
{\bf Proof.} The fact that
\[ [ m , \bar{m} ] \perp {\cal L} \mbox{ and } [ m , \bar{m} ] \perp {\cal N} \]
is a direct consequence of Frobenius's theorem since the planes orthogonal to $\cal L$ and $\cal N$ are integrable. But then we notice that
\[ [ m , \bar{m} ] = \nabla_m \bar{m} - \nabla_{\bar{m}} m = 2 i \Im \left( \nabla_m \bar{m} \right) \, ,\]
so we must have
\[ {\cal L}_a \Im \left( \nabla_m \bar{m}^a \right) = {\cal N}_a \Im \left( \nabla_m \bar{m}^a \right) =0 \]
i.e.
\begin{equation} \label{KillImaginaryPart}
{\cal L}_a \nabla_m \bar{m}^a = {\cal L}_a \nabla_{\bar{m}} m^a \mbox{ and } {\cal N}_a \nabla_m \bar{m}^a = {\cal N}_a \nabla_{\bar{m}} m^a \, .
\end{equation}
Let us now expand the divergence of ${\cal L}$ in the basis ${\cal L}$, $\cal N$, $m$, $\bar{m}$~: 
\begin{eqnarray*}
\mathrm{div} {\cal L} = g_{ab} \nabla^a {\cal L}^b &=& \left( {\cal L}_a {\cal N}_b + {\cal N}_a {\cal L}_b - m_a \bar{m}_b - \bar{m}_a m_b \right) \nabla^a {\cal L}^b \\
&=& {\cal N}_b \nabla_{\cal L} {\cal L}^b + {\cal L}_b \nabla_{\cal N} {\cal L}^b - \bar{m}_b \nabla_{m} {\cal L}^b - m_b \nabla_{\bar{m}} {\cal L}^b \\
&=& - 2 m_b \nabla_{\bar{m}} {\cal L}^b = -2 \rho \, ,
\end{eqnarray*}
the final simplification being obtained since $\nabla_{\cal L} {\cal L} =0$, $\cal L$ is null and by \eqref{KillImaginaryPart} via the following short calculation
\[ \bar{m}_b \nabla_{m} {\cal L}^b + m_b \nabla_{\bar{m}} {\cal L}^b = - {\cal L}^b \left( \nabla_{m} \bar{m}_b + \nabla_{\bar{m}} m_b \right) = -2 {\cal L}^b \nabla_{\bar{m}} m_b = 2 m_b \nabla_{\bar{m}} {\cal L}^b \, .\]
This proves the lemma. \qed
\begin{remark}
Note in particular that $\rho$ is real.
\end{remark}
If we denote by $s$ the affine parameter along the integral curves of $\cal L$, choosing $s=0$ at the vertex, then the behaviour at the vertex of $\rho$ is given by (see Seitz-Schneider-Ehlers \cite{SeSchneEh})\footnote{In \cite{SeSchneEh}, the behaviour obtained is $1/s$ instead of $-1/s$. This change of sign is due to the fact that they chose a vector $\cal L$ which was past pointing whereas our $\cal L$ points to the future. Apart from this sign difference, the normalization conditions for $\cal L$ in \cite{SeSchneEh} are satisfied here at the vertex, in that $g(\nabla u , \nabla t) =1$ there.}~:
\[ \rho = \frac{-1}{s} (1+Ks^2) + O(s^2) \mbox{ as } s \rightarrow 0 \, ,\]
where $K$ is a constant characteristic of the geometry of the vertex.
The affine parameter $s$ can be calculated explicitly in terms of $v$ on the cone~:
\[ \d s (\nabla u ) = 1 \, \mbox{and }
\d v (\nabla u ) = g (\nabla u \, ,~ \nabla v ) = \frac{4}{N^2} \, , \]
hence
\[ \frac{\d s}{\d v} = \frac{N^2}{4} \]
and since $N^2 =2$ at the vertex, on ${\cal C}^+_0$,
\[ s \simeq \frac{v}{2} + O(v^2) =r + O(r^2) \mbox{ as } r \rightarrow 0 \, .\]
So we have
\begin{equation} \label{EquivRho}
\rho = \frac{-1}{r} (1+K r^2) + O(r^2) \mbox{ as } r \rightarrow 0 \, .
\end{equation}
\begin{remark} \label{Smoothnessrrho}
This expansion guarantees enough regularity for our results. However, one would expect that $r \rho$ is in fact s-smooth on $\mathfrak{C}$ and that the expansion given in \cite{SeSchneEh} was stopped at the order $s^2$ to avoid unnecessarily heavy calculations and not because of a fundamental difficulty. Since this in no way affects our results, we shall simply consider that $r \rho$ is s-smooth on $\mathfrak{C}$. The cautious reader is free to consider this smoothness as that given by the expansion at order $r^2$ only.
\end{remark}
We now perform the rescaling of the transverse frame spinor~: we put
\begin{equation} \label{SpinDyadRescaled}
\hat{O}^A := O^A \, ,~ \hat{I}^A := r I^A \, .
\end{equation}
We denote with a hat all quantities referring to the new dyad and its associated tetrad
\[ \hat{\cal L} := {\cal L} \, ,~ \hat{\cal N} := r^2 {\cal N} \, ,~ \hat{m} := r m \, ,\]
and in particular
\[ \hat{\phi_0}:= \phi_0 \, ,~ \hat{\phi}_1 := r \phi_1 \, ,~ \hat{\chi}_{0'} := \chi_{0'} \, ,~\hat{\chi}_{1'} := r \chi_{1'} \, .\]
The constraints are the second and fourth equations in the expression \eqref{NPCompactDirac} of the Dirac system in the compacted spin-coefficient formalism. We multiply them both by $r$. All the terms in the left hand-side of the first equation have weight $\{ 1,1;1,0\}$ and all those in the left hand-side of the second have weight $\{ 1,0;1,1\}$~; under the rescaling \eqref{SpinDyadRescaled}, all these terms thus undergo a scaling by $r$. We can therefore re-interpret the constraints multiplied by $r$ as follows~:
\[
 \begin{array}{l}
 { \hat{\mbox{\th}} \hat\phi_1 - \hat\eth' \hat\phi_0 - \hat\pi \hat\phi_0 - \hat\rho \hat\phi_1 -iq \hat{\cal L}^a \Phi_a \hat\phi_1 + i q \bar{\hat{m}}^a \Phi_a \hat\phi_0 = - \frac{m}{\sqrt{2}} r \hat\chi_{0'} \, , } \\ \\
{ \hat{\mbox{\th}} \hat\chi_{1'} - \hat\eth \hat\chi_{0'} - \bar{\hat\pi} \hat\chi_{0'} - \bar{\hat\rho} \hat\chi_{1'} -iq \hat{\cal L}^a \Phi_a \hat\chi_{1'} + i q \hat{m}^a \Phi_a \hat\chi_{0'} = - \frac{m}{\sqrt{2}} r
\hat\phi_{0} \, .} \end{array}
\]
The first term in the first equation is
\[ \hat{\mbox{\th}} \hat\phi_1 = r {\mbox{\th}} \phi_1 = {\cal L}^a \nabla_a \hat\phi_1 - \hat{\gamma'} \hat{\phi_1} \]
and
\[ \hat{\gamma'} = - \frac{1}{r} \hat{O}^A \nabla_{\cal L} \hat{I}_A = \gamma' + \frac{\nabla_{\cal L} r}{r} \, .\]
Similarly,
\[ \hat{\mbox{\th}} \hat\chi_{1'} = {\cal L}^a \nabla_a \hat\chi_{1'} - \bar{\gamma'} \hat\chi_{1'} - \frac{\nabla_{\cal L} r}{r} \hat\chi_{1'} \, . \]
So the rescaled constraints read
\[
 \begin{array}{l}
 { \left( {\cal L}^a \nabla_a - {\gamma'} \right) \hat\phi_1 - \frac{\nabla_{\cal L} r}{r} \hat\phi_1 - \hat\eth' \hat\phi_0 - \hat\pi \hat\phi_0 - \hat\rho \hat\phi_1 -iq \hat{\cal L}^a \Phi_a \hat\phi_1 + i q \bar{\hat{m}}^a \Phi_a \hat\phi_0 = - \frac{m}{\sqrt{2}} r \hat\chi_{0'} \, , } \\ \\
{ \left( {\cal L}^a \nabla_a - \bar{\gamma'} \right) \hat\chi_{1'} - \frac{\nabla_{\cal L} r}{r} \hat\chi_{1'}- \hat\eth \hat\chi_{0'} - \bar{\hat\pi} \hat\chi_{0'} - \bar{\hat\rho} \hat\chi_{1'} -iq \hat{\cal L}^a \Phi_a \hat\chi_{1'} + i q \hat{m}^a \Phi_a \hat\chi_{0'} = - \frac{m}{\sqrt{2}} r
\hat\phi_{0} \, .} \end{array}
\]
The derivative of $r$ along $\cal L$ can be calculated easily
\[ \nabla_{\cal L} r = \frac12 \d v ({\cal L}) = \frac{2}{N^2} \, , \]
which along a given integral line of $\cal L$ can be expanded as follows near $r=0$~:
\[ \nabla_{\cal L} r = 1 - \frac{4}{N^3(p_0)}  r \nabla_{\cal L} N (p_0) + O(r^2) = 1-\sqrt{2} r \nabla_{\cal L} N (p_0) + O(r^2) \, . \]
The coefficient $\rho$ is of weight $\{ 1,0;1,0\}$ so we have $\hat{\rho} = \rho$ and we see that as $r \rightarrow 0$,
\[ \frac{\nabla_{\cal L} r}{r} + \hat{\rho} \rightarrow -\sqrt{2} \nabla_{\cal L} N (p_0) \]
which depends smoothly on $\omega$. So remembering that $\rho$ is real, we have the new rescaled form of the constraint equations~:
\begin{equation} \label{RescConstraintEq}
\left. \begin{array}{l}
 { ( \frac{4}{N^2} \partial_v - \gamma' ) \hat\phi_1 - \left( \frac{\nabla_{\cal L} r}{r} + \hat{\rho} \right) \hat\phi_1 - ( \hat\eth' + \hat\pi ) \hat\phi_0 -iq \hat{\cal L}^a \Phi_a \hat\phi_1 + i q \bar{\hat{m}}^a \Phi_a \hat\phi_0 = - \frac{m}{\sqrt{2}} r \hat\chi_{0'} \, , } \\ \\
{ ( \frac{4}{N^2} \partial_v -\bar{\gamma'} \, ) \hat{\chi}_{1'} - \left( \frac{\nabla_{\cal L} r}{r} + \hat{\rho} \right) \hat\chi_{1'}- ( \hat\eth + \bar{\hat\pi} ) \hat\chi_{0'} -iq \hat{\cal L}^a \Phi_a \hat\chi_{1'} + i q \hat{m}^a \Phi_a \hat\chi_{0'} = - \frac{m}{\sqrt{2}} r
\hat\phi_{0} \, .} \end{array} \right\}
\end{equation}
\begin{lemma} \label{epsilonpi}
The coefficients $\gamma'$ and $\pi$ are s-smooth on $\mathfrak{C}$ and $\hat{\pi} \rightarrow 0$ at the vertex.
\end{lemma}
{\bf Proof.} Recall that the tetrad $\{ l,n,m,\bar{m} \}$ is s-smooth on $\mathfrak{C}$ and therefore so are the tetrads $\{ {\cal L} , {\cal N}, m , \bar{m} \}$ and $\{ \hat{\cal L} , \hat{\cal N}, \hat{m} , \bar{\hat{m}} \}$. In particular any directional derivative along $l$ or $\hat{\cal L} = {\cal L}$ of any of the frame vectors above will be s-smooth on $\mathfrak{C}$. The expressions of $\gamma'$ and $\pi$ are
\begin{gather*}
\gamma' = - \varepsilon = - \frac{1}{2} \left( {\cal N}^a \nabla_{\cal L} {\cal L}_a + m^a \nabla_{\cal L} \bar{m}_a \right) = -\frac12 m^a \nabla_{\cal L} \bar{m}_a \, , \\
\pi = I^A \nabla_{\cal L} I_A = - {\bar{m}}^a \nabla_{\cal L} {\cal N}_a = {\cal N}^a \nabla_{\cal L} {\bar{m}}_a \, .
\end{gather*}
They only involve derivatives along $\cal L$ of the frame vector $\bar{m}^a$ and are therefore s-smooth on $\mathfrak{C}$ (see lemma \ref{Sregl}). Moreover, we have
\[ \hat\pi = \frac{1}{r} \hat{I}^A \nabla_{\cal L} \hat{I}_A = r \pi + I^A I_A \nabla_{\cal L} r = r \pi \, .\]
This concludes the proof. \qed

The angular operators $\hat{\eth}'$ and $\hat{\eth}$, acting respectively on $\hat{\phi}_0$ and $\hat{\chi}_{0'}$, involve covariant derivatives along $\bar{\hat{m}}$ and $\hat{m}$ as well as some spin-coefficients~:
\[ \hat{\eth}' \hat{\phi}_0 = \left( \nabla_{\bar{\hat{m}}} - \hat{\alpha} \right) \hat{\phi}_0 \, ,~ \hat{\eth} \hat{\chi}_{0'} = \left( \nabla_{\hat{m}} - \bar{\hat{\alpha}} \right) \hat{\chi}_{0'} \, .  \]
The components $\hat{\phi}_0 = \phi_0 = \phi_A O^A$ and $\hat{\chi}_{0'} = \chi_{0'} = \chi_{A'} \bar{O}^{A'}$ are s-smooth functions on $\mathfrak{C}$, i.e. smooth functions of $(r,\omega)$ on $[0,T] \times S^2$. Corollary \ref{ConeSmoothCylinder} then implies that $\nabla_{\bar{\hat{m}}} \hat{\phi}_0$ and $\nabla_{\hat{m}} \hat{\chi}_{0'} $ are also s-smooth functions on $\mathfrak{C}$ (once again ignoring the necessity to work with two charts). We then need to deal with the coefficient $\hat{\alpha}$. However, keeping the whole quantities $\hat{\eth}' \hat{\phi}_0$ and $\hat{\eth} \hat{\chi}_{0'}$ makes the proofs easier and as an upshot allows to establish the regularity of $\hat{\alpha}$, and more precisely of both terms involved in $\hat{\alpha}$, on $\mathfrak{C}$ which will be useful for appendix \ref{MoreSpinCoeffs}.

Recall that there exists a smooth spinor field $F_A \oplus G^{A'}$ on $D_T$ such that on ${\cal C}^+_0$ we have
\[ F_0 = \phi_0 \mbox{ and } G_{0'} = \chi_{0'} \, .\]
\begin{lemma} \label{RegularityAngularDerivatives}
The quantities $\hat{\eth}' \hat{\phi}_0$ and $\hat{\eth} \hat{\chi}_{0'}$ are s-smooth on $\mathfrak{C}$.
\end{lemma}
{\bf Proof.} We write the proof for $\hat{\eth}' \hat{\phi}_0$~: it can be considered as the action of $r\eth'$ on $F_0$ and this can be treated as follows (the passage from the third to the fourth line of the calculation below uses the expression of the covariant derivative along $\bar{m}$ of the first basis spinor given in Penrose-Rindler \cite{PeRi} Vol. 1, p. 227, equation (4.5.26))~:
\begin{eqnarray}
\hat{\eth}' \hat{\phi}_0 &=& r \eth' \phi_0 = r \eth' F_0 \nonumber \\
&=& r \left( \nabla_{\bar{m}} - \alpha \right) F_0= r \nabla_{\bar{m}} \left( F_A O^A \right) - r \alpha F_0 \nonumber \\
&=& r \left( \nabla_{\bar{m}} F_A \right) O^A + r F_A \nabla_{\bar{m}} O^A - r \alpha F_0 \nonumber \\
&=& r \left( \nabla_{\bar{m}} F_A \right) O^A + r F_A \left( \alpha O^A - \rho I^A \right) - r \alpha F_0 \nonumber \\
&=& r \left( \nabla_{\bar{m}} F_A \right) O^A - r \rho F_1 \, . \label{CalculEth'phi0}
\end{eqnarray}
This is s-smooth on $\mathfrak{C}$ by regularity of $F_A$, $\bar{m}$, $O^A$, $I^A$ and $r\rho$ (see remark \ref{Smoothnessrrho}). The treatment for $\hat{\eth} \hat{\chi}_{0'}$ is similar. \qed

As mentioned above, we now give a more precise regularity result for $\hat{\alpha}$~; this is not necessary for the present proof but will be useful in appendix \ref{MoreSpinCoeffs}.
\begin{corollary} \label{TermsAlphaHat}
The two terms $r^{-2} \hat{\cal N}^a \nabla_{\bar{\hat{m}}} {\cal L}_a$ and $r^{-2}\hat{m}^a \nabla_{\bar{\hat{m}}} \bar{\hat{m}}_a$ of which $\hat{\alpha}$ is composed\footnote{The complete expression of $\hat{\alpha}$ is
\[ \hat{\alpha} = \frac{1}{2r^2} \left( \hat{\cal N}^a \nabla_{\bar{\hat{m}}} {\cal L}_a + \hat{m}^a \nabla_{\bar{\hat{m}}} \bar{\hat{m}}_a + r \nabla_{\bar{\hat{m}}} r \right) \]
but $\nabla_{\bar{\hat{m}}} r =0$.} are s-smooth on $\mathfrak{C}$.
\end{corollary}
{\bf Proof.} First, note that $\hat{\alpha} \hat{\phi}_0$ is s-smooth of $\mathfrak{C}$ as a consequence of lemma \ref{RegularityAngularDerivatives}, of the equality
\[ \hat{\eth}' \hat{\phi}_0 = \left( \nabla_{\bar{\hat{m}}} - \hat{\alpha} \right) \hat{\phi}_0 \]
and of the s-regularity of $\nabla_{\bar{\hat{m}}} \hat{\phi}_0$ on $\mathfrak{C}$ (which follows from that of $\hat{\phi}_0$ and of $\bar{\hat{m}}$). This is for any given data $\hat{\phi_0} = \phi_0$ obtained from a smooth spinor field, so $\hat{\alpha}$ is s-smooth on $\mathfrak{C}$. Hence we only need to establish the regularity of one of the two terms. We work with the second~: $r^{-2}\hat{m}^a \nabla_{\bar{\hat{m}}} \bar{\hat{m}}_a$. Due to the scalar product with $\hat{m}^a$, only the part of $\nabla_{\bar{\hat{m}}} \bar{\hat{m}}_a$ tangent to the sphere $S_r$ is involved, so
\[ r^{-2}\hat{m}^a \nabla_{\bar{\hat{m}}} \bar{\hat{m}}_a = r^{-2}\hat{m}^a D_{\bar{\hat{m}}} \bar{\hat{m}}_a \]
where $D$ is the Levi-Civita connection on $( S_r , -g |_{S_r} )$. But then $\mathfrak{g} |_{S_r} = -r^2 g |_{S_r}$ and since $r$ is constant on $S_r$, this does not change the covariant derivative. It follows that, denoting by $\mathfrak{D}$ the Levi-Civita connection on $( \mathfrak{C} , \mathfrak{g} )$,
\[ r^{-2}\hat{m}^a \nabla_{\bar{\hat{m}}} \bar{\hat{m}}_a = r^{-2}\hat{m}^a \mathfrak{D}_{\bar{\hat{m}}} \bar{\hat{m}}_a = \mathfrak{g}_{ab} \hat{m}^a \mathfrak{D}_{\bar{\hat{m}}} \bar{\hat{m}}^b \]
and this, by corollary \ref{ConeSmoothCylinder}, is a smooth function of $(r,\omega )$ on $[0,T] \times S^2 $. This concludes the proof. \qed

It follows that the rescaled constraint equations on the cone can be integrated from any data for $\hat{\phi_1}$ and $\hat{\chi}_{1'}$ at the vertex. Of course the only data compatible with the continuity of the physical field at the vertex are zero for both $\hat{\phi_1}$ and $\hat{\chi}_{1'}$. It remains to understand the behaviour of the components ${\phi_1}$ and ${\chi}_{1'}$ of the physical field at the vertex. The value of $\phi_1$ at the vertex (at least along a given future null geodesic) can be recovered from the value of $\partial_v \hat{\phi}_1$ at $p_0$ along the same geodesic, indeed
\[ \lim_{r\rightarrow 0} {\phi}_1 = \lim_{r\rightarrow 0} \frac{\hat{\phi}_1}{r} = 2\lim_{v\rightarrow 0} \frac{\hat{\phi}_1}{v} = 2 \partial_v \hat{\phi}_1 (p_0) \, .\]
This value can be extracted from the constraint equation~:
\[ ( \frac{4}{N^2} \partial_v - \gamma' ) \hat\phi_1 - \left( \frac{\nabla_{\cal L} r}{r} + \hat{\rho} \right) \hat\phi_1 - ( \hat\eth' + \hat\pi ) \hat\phi_0 -iq \hat{\cal L}^a \Phi_a \hat\phi_1 + i q \bar{\hat{m}}^a \Phi_a \hat\phi_0 = - \frac{m}{\sqrt{2}} r \hat\chi_{0'} \, .\]
We can infer the behaviour of each term as $r\rightarrow 0$ from the results we have already proven. First, $\gamma'$ is s-smooth on $\mathfrak{C}$ and $\hat{\phi}_1$ tends to zero at the vertex, so $\gamma' \hat{\phi}_1$ tends to zero as $r\rightarrow 0$. The same is true for the term $\left( \frac{\nabla_{\cal L} r}{r} + \hat{\rho} \right) \hat\phi_1$. Similarly $\hat{\pi}$ tends to zero and $\hat{\phi}_0 = \phi_0$ has a finite limit at the vertex (along each future null geodesic), so $\hat{\pi} \hat{\phi}_0 $ tends to zero at the vertex. Then ${\cal L}^a$ is a bounded vector field and $\Phi_a$ a smooth $1$-form on $D_T$, so the function ${\cal L}^a \Phi_a$ is bounded on $D_T$ and when we multiply it with $\hat{\phi}_1$, the product tends to zero at the vertex. The vector field $m^a$ is bounded, so $\hat{m}^a \Phi_a \hat{\phi}_0 = r m^a \Phi_a \phi_0$ also tends to zero at the vertex. Finally $r \chi_{0'}$ tends to zero as $r \rightarrow 0$. So we see that the value at the vertex reached along a given future null geodesic of $\phi_1$ is given by (remember that $4/N^2$ is equal to $2$ at $p_0$)
\[ \phi_1 (0,\omega ) = 2 \partial_v \hat{\phi}_1 (0,\omega ) = \lim_{r \rightarrow 0} \hat{\eth}' \hat{\phi}_0 (r,\omega) \, .\]
Now, recall the expression \eqref{CalculEth'phi0} of $\hat{\eth}' \hat{\phi}_0$ obtained in the proof of lemma \ref{RegularityAngularDerivatives}~:
\[ \hat{\eth}' \hat{\phi}_0 = r \left( \nabla_{\bar{m}} F_A \right) O^A - r \rho F_1 \, .\]
The first term on the right-hand side tends to zero as $r \rightarrow 0$ by regularity of $F_A$, $\bar{m}$ and $O^A$, so
\[ \phi_1 (0,\omega ) = \lim_{r\rightarrow 0} \hat{\eth}' \hat{\phi}_0 (r,\omega) = -\lim_{r\rightarrow 0} r \rho F_1(r,\omega) = F_1 (0, \omega ) \, .\]
With a similar proof we obtain
\[ \chi_{1'} (0,\omega ) = G_{1'} (0,\omega ) \, .\]
In the spin-frame $\{ o^A \, ,~ \iota^A \}$, which differs from $\{ O^A \, ,~ I^A \}$ only by a smooth scaling, we still have at the vertex
\[ \phi_1 (0,\omega ) = F_1 (0, \omega ) \, ,~\chi_{1'} (0,\omega ) = G_{1'} (0,\omega ) \, .\]
Since the Dirac spinor $F_A \oplus G^{A'}$ is continuous at $p_0$, it follows that $\phi_1$ and $\phi_0$ and also $\chi_{1'}$ and $\chi_{0'}$, all being evaluated in the spin-frame $\{ o^A \, ,~\iota^A \}$, must be related by the continuity matching conditions at the vertex. This concludes the proof of proposition \ref{WellPosedConstraintEq}. \qed

\section{Behaviour of spin coefficients} \label{MoreSpinCoeffs}

In this appendix, we study the behaviour of the spin coefficients involved in the definition of ${\cal H}_{{\cal C}^+_0}$ (section \ref{H1Setting}, definition \ref{H1DataSpace}) that have not yet been studied (namely $\alpha$, $\beta$, $\kappa$, $\mu$ and $\sigma$ in the Newman-Penrose tetrad $\{ l , n , m , \bar{m} \}$), and also the coefficient $\nu$ which appears in section \ref{secProofT2}. Then we give a summary of the behaviour at the vertex of all the spin coefficients that we use.

We start with $\mu$ and $\beta$. We work out their behaviour using the tetrad $\{ {\cal L} , {\cal N} , m , \bar{m} \}$. The two spin-frames differ only by a scaling of $l$ and $n$ by a smooth function on $D_T$, so the behaviour of the coefficients worked out in one tetrad will readily give the behaviour in the other. The expressions of these two coefficients are
\begin{eqnarray*}
\mu &=& I^A \nabla_m I_A = - \bar{m}^a \nabla_m {\cal N}_a = {\cal N}^a \nabla_m \bar{m}_a = - \rho' \, , \\
\beta &=& I^A \nabla_m O_A = \frac12 \left( {\cal N}^a \nabla_m {\cal L}_a + m^a \nabla_m \bar{m}_a \right) \, .
\end{eqnarray*}
The behaviour of $\rho'$ was studied in the proof of proposition \ref{GaussCurvature}, it is equivalent to $1/r$ at the tip of the cone. The essential argument was that swapping $l$ and $n$ corresponds to exchanging the future and past null cones together with a regular dependence on the tip of the cone.  Note that the same argument can be used to study the behaviour of the two coefficients $\gamma$ and $\tau$ (which appear in the expression $(\nabla_n \Psi)_{1,4}$ in section \ref{H1Setting} but are cancelled in the expression of the norm by terms from the equation). They are obtained from $\varepsilon$ and $\pi$ by swapping $l$ and $n$ and changing the sign. We have already established the s-regularity on $\mathfrak{C}$ of $\varepsilon = -\gamma'$ and $\pi$, the same is therefore true for $\gamma$ and $\tau$. We cannot use the same argument for $\beta$ in relation to $\alpha$ but the calculations we have done for $\hat{\alpha}$ will be useful here nonetheless. The rescaled coefficient $\hat{\beta}$ is
\begin{eqnarray*}
\hat{\beta} &=& \frac{1}{2r^2} \left( \hat{\cal N}^a \nabla_{\hat{m}} \hat{\cal L}_a + \hat{m}^a \nabla_{\hat{m}} \bar{\hat{m}}_a + r \nabla_{\hat{m}} r \right) \\
&=& \frac{1}{2r^2} \left( \hat{\cal N}^a \nabla_{\hat{m}} \hat{\cal L}_a + \hat{m}^a \nabla_{\hat{m}} \bar{\hat{m}}_a \right) = r \beta \, .
\end{eqnarray*}
The same reasoning as in the proof of corollary \ref{TermsAlphaHat} gives the s-regularity of $r^{-2} \hat{m}^a \nabla_{\hat{m}} \bar{\hat{m}}_a$ on $\mathfrak{C}$. Another conclusion of corollary \ref{TermsAlphaHat} was the s-regularity of $r^{-2} \hat{\cal N}^a \nabla_{\bar{\hat{m}}} \hat{\cal L}_a$ which is the complex conjugate of the first term of $\beta$. This proves the s-regularity of $\hat{\beta} = r \beta$ on $\mathfrak{C}$.

The coefficient $r\alpha = \hat{\alpha}$ is smooth near the vertex by corollary \ref{TermsAlphaHat}. From \cite{SeSchneEh}, the coefficient $\sigma /r$ is smooth near the vertex. As for $\kappa = m^a \nabla_{\cal L} {\cal L}_a$, it is zero because $\nabla_{\cal L} {\cal L} =0$. Note that in the tetrad $\{ l,n,m,\bar{m} \}$, $\nabla_l l$ is proportional to $l$ and therefore orthogonal to $m$, so we still have $\kappa = m^a \nabla_l l_a =0$. The same is true of $\nu = - \kappa' = -\bar{m}^a \nabla_{\cal N} {\cal N}_a =0$ since $\nabla_{\cal N} {\cal N}$ is parallel to $\cal N$ and therefore orthogonal to $\bar{m}$ (and also in the tetrad $\{ l,n,m,\bar{m} \}$).

We now summarize the behaviour at the vertex of all the spin coefficients we use, defined in the tetrad $\{ {\cal L} , {\cal N} , m , \bar{m} \}$.
\begin{lemma} \label{SpinCoeffs}
The functions $\frac{1}{r} \sigma $, $\varepsilon$, $\gamma$, $\pi$, $\tau$, $r \alpha$, $r\beta$ and $r\mu$ are smooth on the cone, $\kappa = \nu =0$ everywhere and $\rho = -\frac{1}{r} (1+Kr^2) +O(r^2 )$ near the vertex.
\end{lemma}

\section{Compacted spin coefficient formalism} \label{CompactedNP}

The behaviour of scalars under a rescaling of the frame spinors depends on their weight, which is a collection of $4$ integers (or possibly real numbers)~: a scalar $\eta$ is said to have weight $\{ r',r;t',t \}$ if under a rescaling of the spin-frame by nowhere vanishing scalar fields $\lambda$ and $\mu$,
\[ o^A \mapsto \lambda o^A \, ,~ \iota^A \mapsto \mu \iota^A \, ,\]
it transforms as
\[ \eta \mapsto \lambda^{r'} \mu^r \bar{\lambda}^{t'} \bar{\mu}^t \eta \, .\]
For example, $\phi_0 = \phi_A o^A$ has weight $\{ 1,0;0,0\}$ and $\phi_1 = \phi_A \iota^A$ has weight $\{ 0, 1 ; 0 , 0 \}$. If we work with normalized spin-frames, then to preserve the normalization we must impose $\mu = 1/\lambda$ and then only two numbers play a role~: $p = r'-r$ and $q=t'-t$. A scalar is then said to have weight $\{ p;q \}$ or equivalently to have boost weight $\frac12 (p+q)$ and spin weight $\frac12 (p-q)$.
Not all scalars have a weight and the fundamental directional derivatives of the Newman-Penrose formalism $l^a \partial_ a$, $n^a \partial_ a$, $m^a \partial_ a$ and $\bar{m}^a \partial_a$ do not transform weighted scalars into weighted scalars. The compacted spin-coefficient formalism groups these derivatives with unweighted spin-coefficients to produce weighted operators denoted $\thorn$, $\thorn '$, $\eth$, $\eth '$ (pronounced ``thorn'' and ``eth''). The action of these operators on a weighted scalar $\eta$ of weight $\{ r',r;t',t \}$ is defined by
\begin{eqnarray*}
\mbox{\th} \eta &:=& \left( l^a \partial_a - r' \varepsilon - r \gamma ' - t ' \bar{\varepsilon} - t \bar{\gamma}' \right) \eta \, ,\\
\eth \eta &:=& \left( m^a \partial_a - r' \beta - r \alpha ' - t' \bar{\alpha} - t \bar{\beta}' \right) \eta \, , \\
\eth' \eta &:=& \left( \bar{m}^a \partial_a - r' \alpha - r \beta' - t'\bar{\beta} - t\bar{\alpha}' \right) \eta \, ,\\
\mbox{\th}' \eta &:=& \left( n^a \partial_a - r' \gamma - r \varepsilon'- t'\bar{\gamma} - t \bar{\varepsilon}' \right) \eta \, .
\end{eqnarray*}
The results are weighted scalars with the following respective weights~: $\mbox{\th} \eta$ has weight $\{ r'+1,r;t'+1,t \}$, $\mbox{\th}' \eta$ has weight $\{ r',r+1;t',t+1 \}$, $\eth \eta$ has weight $\{ r'+1,r;t',t+1 \}$ and $\eth' \eta$ has weight $\{ r',r+1;t'+1,t \}$. In the normalized case, we have the following equalities between spin-coefficients~:
\begin{gather*}
\kappa = -\nu' \, ,~ \rho = - \mu' \, ,~ \sigma = -\lambda ' \, ,~ \tau = - \pi ' \, ,~ \varepsilon = - \gamma ' \, ,~ \alpha = - \beta' \, , \\
\kappa ' = -\nu \, ,~ \rho ' = - \mu \, ,~ \sigma ' = -\lambda  \, ,~ \tau '= - \pi  \, ,~ \varepsilon '= - \gamma  \, ,~ \alpha '= - \beta \, 
\end{gather*}
and the weighted derivatives take the simplified expression
\begin{eqnarray*}
\mbox{\th} \eta &:=& \left( l^a \partial_a + p \gamma ' +q \bar{\gamma}' \right) \eta \, ,\\
\eth \eta &:=& \left( m^a \partial_a - p \beta + q \bar{\beta}' \right) \eta \, , \\
\eth' \eta &:=& \left( \bar{m}^a \partial_a +p  \beta' - q\bar{\beta} \right) \eta \, ,\\
\mbox{\th}' \eta &:=& \left( n^a \partial_a - p \gamma - q\bar{\gamma} \right) \eta \, .
\end{eqnarray*}

The expression \eqref{NewmanDirac} of the Dirac equation in the Newman-Penrose formalism is only valid for a normalized spin-frame, so we must use the simplified expressions of the weighted derivatives to obtain the corresponding expression in the compacted spin coefficient formalism. The components of the Dirac field are weighted scalars of weights $\{ 1 ; 0\}$ for $\phi_0$, $\{ -1;0\}$ for $\phi_1$, $\{ 0;1\}$ for $\chi_{0'}$ and $\{ 0;-1\}$ for $\chi_{1'}$. So \eqref{NewmanDirac} can be reformulated as
\begin{eqnarray}
\label{NPCompactDirac}
 \left. \begin{array}{l}
{ \mbox{\th}' \phi_0 - \eth \phi_1 + \mu \phi_0 + \tau \phi_1 - i q n^a \Phi_a \phi_0 + iq m^a \Phi_a \phi_1 = \frac{m}{\sqrt{2}} \chi_{1'} \, , } \\ \\
{ \mbox{\th} \phi_1 - \eth' \phi_0 - \pi \phi_0 - \rho \phi_1 -iq l^a \Phi_a \phi_1 + i q \bar{m}^a \Phi_a \phi_0 = - \frac{m}{\sqrt{2}} \chi_{0'} \, , } \\ \\
{ \mbox{\th}' \chi_{0'} - \eth' \chi_{1'} + \bar{\mu} \chi_{0'}
+ \bar{\tau} \chi_{1'} -iq n^a \Phi_a \chi_{0'} + iq \bar{m}^a \Phi_a \chi_{1'}= \frac{m}{\sqrt{2}} \phi_{1} \,
, } \\ \\
{ \mbox{\th} \chi_{1'} - \eth \chi_{0'} - \bar{\pi} \chi_{0'} - \bar{\rho} \chi_{1'} -iq l^a \Phi_a \chi_{1'} + i q m^a \Phi_a \chi_{0'} = - \frac{m}{\sqrt{2}}
\phi_{0} \, .} \end{array} \right\}
\end{eqnarray}
The second and fourth equations are the constraint equations along the cone.

\section*{Acknowledgments}

Both authors would like to thank the Institut Mittag-Leffler and the organizers of the ``Geometry, analysis and general relativity'' semester during which we started working on this paper. This work was partially supported by the ANR projects JC0546063 and 08-BLAN-0228. We would also like to thank J\'er\'emie Joudioux for fruitful discussions. We are grateful to the referee for pointing out an error in a first version of the paper.


\begin{thebibliography}{100}


\bibitem{BeSa} A.N. Bernal, M. Sanchez, {\em On smooth Cauchy hypersurfaces and Geroch's splitting Theorem}, Comm. Math. Phys., {\bf 243} (2003), 461--470.

\bibitem{Ch}
S. Chandrasekhar, The mathematical theory of black holes. International Series of Monographs on Physics, 69. Oxford Science Publications. The Clarendon Press, Oxford University Press, New York, 1983. 

\bibitem{Fri75} F.G. Friedlander, The wave equation on a curved spacetime, Cambridge Monographs on Mathematical Physics, Cambridge University Press, 1975.

\bibitem{Fri} T. Friedrich, Dirac operators in Riemannian geometry. Translated from the 1997 German original by Andreas Nestke. Graduate Studies in Mathematics, 25. American Mathematical Society, Providence, RI, 2000.

\bibitem{Ge1} R.P. Geroch, {\em Spinor structure of space-times
in general relativity I}, J. Math. Phys. {\bf 9} (1968), 1739--1744.

\bibitem{Ge2} R.P. Geroch, {\em Spinor structure of space-times
in general relativity II}, J. Math. Phys. {\bf 11} (1970), 342--348.

\bibitem{Ge3} R.P. Geroch, {\em The domain of dependence},
      J. Math. Phys. {\bf 11} (1970), 437--449.

\bibitem{Ha} D. H\"afner, {\em Creation of fermions by rotating charged black holes}, M\'emoires de la SMF {\bf 117} (2009), 158 pp.

\bibitem{HoBook} L. H\"ormander, The analysis of linear partial differential operators Vol. I, Springer, 1983.

\bibitem{Ho} L. H\"ormander, {\em A remark on the characteristic Cauchy problem},  J. Funct. Anal.  {\bf 93}  (1990),  no. 2, 270--277.

\bibitem{Jo} J. Joudioux, {\em Integral formula for the characteristic Cauchy problem on a curved background}, arXiv:0910.4620, to appear in J. Math. Pures et Appliquées.

\bibitem{KlaNi} S. Klainerman, F. Nicol\`o, The evolution problem in general relativity
Progress in Mathematical Physics 25, Birkhäuser Boston Inc., Boston, MA, 2003.

\bibitem{MaNi} L.J. Mason, J.-P. Nicolas, {\em Conformal scattering and the Goursat problem},  J. Hyperbolic Differ. Equ. {\bf 1} (2004),  no. 2, 197--233.

\bibitem{Ni} J.-P. Nicolas, {\em Dirac fields on asymptotically flat spacetimes}, Dissertationes Math. (Roz\-prawy Mat.) 408 (2002), 85 pp.

\bibitem{Pe63} R. Penrose, {\em Null Hypersurface Initial Data for Classical Fields of Arbitrary Spin and for General Relativity},  Gen. Relativity Gravitation {\bf 12} (1980), no. 3, 225--264, originally appeared in 1963 in Aerospace Research Laboratories Tech. Documentary Report 63--56 (ed. P. G. Bergmann).

\bibitem{PeRi} R. Penrose and W. Rindler, Spinors and space-time, Vol. I (1984) and Vol. 2 (1986), Cambridge University Press.

\bibitem{Re} A. Rendall, {\em Reduction of the characteristic initial value problem to the Cauchy problem and its applications to the Einstein equations}, Proc. Roy. Soc. London Series A Mathematical and Physical sciences {\bf 427} (1990), no. 1872, 221--239.

\bibitem{SeSchneEh} S. Seitz, P. Schneider, J. Ehlers, {\em Light propagation in arbitrary spacetimes and the 
gravitational lens approximation}, Class. Quantum Grav. {\bf 11} (1994), 2345-2373.

\bibitem{Stie} E. Stiefel, {\em Richtungsfelder und Fernparallelismus in $n$-dimensionalen Mannigfaltigkeiten}, Comment. Math. Helv. {\bf 8} (1936), 305--353.

\bibitem{Wa} R. Wald, General relativity, The University of Chicago Press, 1984.

\end{thebibliography}
\end{document}